%% file: main.tex
\pdfoutput=1
\documentclass[11pt]{article}
\usepackage[utf8]{inputenc}

\usepackage[margin=1in]{geometry}
\usepackage{authblk}

\usepackage{tikz}
\usetikzlibrary{arrows.meta, positioning}

\usepackage{graphicx}
\usepackage{enumerate,bbm}
\usepackage[numbers,sort&compress]{natbib}
\usepackage{xcolor,hyperref}
\usepackage{microtype}
\usepackage{amsthm,amsfonts,amssymb,amsmath}
\usepackage{booktabs,array}
\usepackage[capitalize]{cleveref}
\usepackage[skip=0pt]{caption}

\hypersetup{colorlinks,allcolors=blue,linktocpage=true}
\theoremstyle{definition}
\newtheorem{definition}{Definition}
\newtheorem{lemma}{Lemma}
\newtheorem{proposition}{Proposition}
\newtheorem{corollary}{Corollary}
\newtheorem{remark}{Remark}
\newtheorem{assumption}{Assumption}
\newtheorem{theorem}{Theorem}

\renewcommand{\S}{\mathcal S}
\newcommand{\N}{{\cal N}}
\newcommand{\ho}{\mathbb}
\newcommand{\ind}{\perp\!\!\!\perp}
\newcommand{\rb}{\right}
\newcommand{\lb}{\left}
\newcommand{\wh}{\widehat}
\newcommand{\bdefn}{\begin{definition}}
\newcommand{\edefn}{\end{definition}}
\newcommand{\bcoro}{\begin{corollary}}
\newcommand{\ecoro}{\end{corollary}}
\newcommand{\sse}{\subseteq}

\begin{document}

\title{Clustered Switchback Designs for Experimentation Under Spatio-temporal Interference}

\author[1]{Su Jia\thanks{Corresponding author.}}
\author[2]{Nathan Kallus}
\author[2]{Christina Lee Yu}

\affil[1]{Center for Data Science for Enterprise and Society (CDSES), Cornell University\protect\\ \texttt{sj693@cornell.edu}}
\affil[2]{School of Operations Research and Information Engineering (ORIE), Cornell University\protect\\ \texttt{\{kallus, cleeyu\}@cornell.edu}}

\date{}

\maketitle

\begin{abstract}
Consider experimentation in the presence of non-stationary dynamics and outcomes, spillover effects (spatial interference), and carryover effects (temporal interference). Given the outcomes measured from the experiment, we want to estimate the global average treatment effect (GATE), i.e. the difference between average outcomes when all units for the entire time horizon are jointly exposed to either treatment or control. We model spatial interference by a graph, where a unit's outcome depends on its neighborhood's treatments, and we model temporal interference by a Markov decision process (MDP), where the transition kernels under either treatment or control arms satisfy a rapid mixing condition. While cluster randomized designs are commonly used to address spillover effects, and switchback designs are commonly used to address carryover effects, we propose a clustered switchback design, where units are grouped into clusters and time steps are grouped into blocks, and each cluster-block combination is jointly assigned to treatment or control.
Under this design, we show that a truncated Horvitz-Thompson estimator performs well, as long as the graph admits a good clustering or the Markovian transitions mix quickly.
Our results provide guidance for practitioners to jointly select good spatial clusters and temporal blocks that yield good performance. We generalize results in the literature that consider either only spatial interference or temporal interference in isolation. By jointly exploiting both spatial and temporal interference we achieve improved performance, where the MSE scales as the product of the MSE from randomizing only over space or time. Simulation studies validate the favorable performance and practical value of our approach.
\end{abstract}

\noindent\textbf{Keywords:} Markov Decision Process (MDP), Average Treatment Effect (ATE), Horvitz-Thompson Estimator, Clustered Randomization

\noindent\textbf{MSC:} 62L05, 62K10, 62M02, 62B15

\section{Introduction}\label{sec:intro}
\input{sections/intro}

\section{Model Setup and Experimental Design}
\label{sec:model}
\input{sections/model}

\section{Main Result}\label{sec:results}
\input{sections/main_results}

\section{Simulation Study}
\label{label:simulation}
\input{sections/xp}

\bibliography{references}
\bibliographystyle{abbrvnat}
\appendix
\input{sections/apdx}
\end{document}

%% file: sections/intro.tex
Randomized experimentation, or A/B testing, is widely used to estimate causal effects on online platforms.
Basic strategies involve partitioning the experimental units (e.g., individuals or time periods) into two groups randomly, and assigning one to treatment and the other to control.
A key challenge is interference: From two-sided markets to social networks, interference between individuals complicates experimentation and makes it difficult to estimate the true effect of a treatment.

Spillover effects in experimentation have been extensively studied \citep{manski2013identification, aronow2017estimating, li2021causal,ugander2013graph, SussmanAiroldi17, toulis2013estimation, BasseAiroldi15, cai2015social, GuiXuBhasinHan15,EcklesKarrerUgander17, chin2019regression}.
Most of these assume neighborhood interference, where the spillover effect is constrained to the direct neighborhood of an individual as given by an {\em interference graph}.
Under this assumption, \citet{ugander2013graph} proposed a clustering-based design, and showed that if the growth rate $\kappa$ of neighborhoods is bounded, then the Horvitz-Thompson (HT) estimator achieves a mean squared error (MSE) of $O(d) \cdot 2^{O(\kappa)}/N$ where $d$ is the maximum degree.
\citet{ugander2023randomized} showed that by introducing randomness into the clustering, the dependence on $\kappa$ can be improved to polynomial.
As there are many settings in which interference extends beyond direct neighbors, \citet{leung2022rate} considers a relaxed assumption in which the interference is not restricted to direct neighbors, but decays as a function of the spatial distance between individuals with respect to an embedding of the individuals in Euclidean space.

Orthogonal to spillover effect, the {\em carryover} effect (or {\em temporal} interference), where past treatments may affect future outcomes, has also been extensively studied.
\citet{bojinov2023design} considered a model in which the temporal interference is bounded by a fixed window length.
Others modeled temporal interference using Markovian evolution of {\em states}, which allows for interference effects that can persist across long time horizons \citep{glynn2020adaptive,farias2022markovian,hu2022switchback,johari2022experimental,shi2023dynamic}.
A common design in practice is \emph{switchback experiments}: The exposure of the entire system (viewed as a single experimental unit) alternates randomly between treatment and control for sufficiently long contiguous blocks of time such that the temporal interference around the switching points does not dominate. \citet{hu2022switchback} showed that a $\tilde O(t_{\rm mix}/T)$ MSE can be achieved where $t_{\rm mix}$ is  (loosely) the mixing time of the Markov chains.

While prior studies either focused on only network interference or only temporal interference, there are many practical settings in which both types of interference are present, such as online platforms, healthcare systems, or ride-sharing networks.
An individual’s outcome may depend not only on who else is treated nearby but also on how the individual's ``state'' has evolved over time, making it essential to develop methodologies that can handle both dimensions jointly.

\begin{figure}
    \centering    \includegraphics[width=0.6\linewidth]{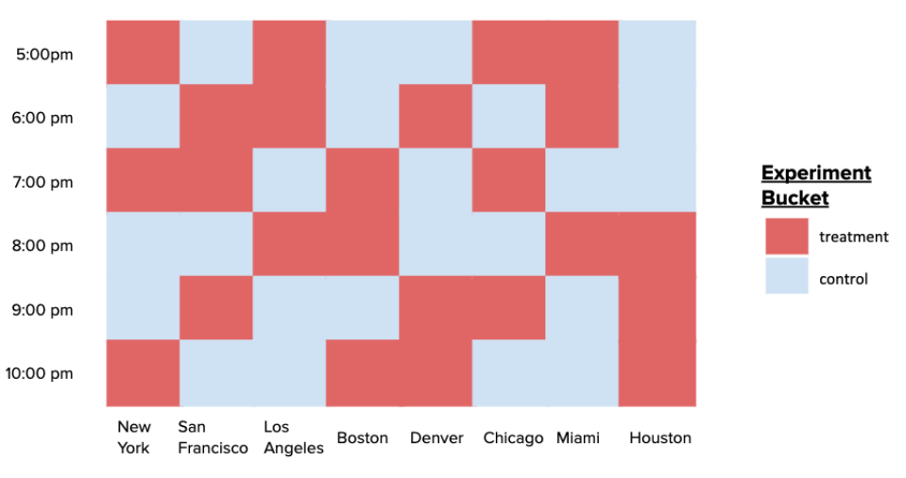}
    \caption{{\bf Clustered Switchback Experiments}.
    The image illustrates clustered switchback on DoorDash \citep{sneider19}.
    The time and geographical locations are grouped into blocks.
    Each spatio-temporal cluster (i.e., Cartesian product) is independently assigned a treatment arm.
    The goal is to estimate the difference in the average (counterfactual) ``outcomes'' (e.g., revenues) between the all-treatment and all-control policy.}
    \label{fig:CSB}
\end{figure}

To address this, {\em clustered switchback experiments} have been recently explored in industry.
The idea is to partition both the space (e.g., by geographics or a social network) and time into discrete blocks, and then randomize at the level of space-time {\em clusters} (i.e., product set of spatial and temporal blocks).
For example, DoorDash randomizes promotions at the region-hour level (see figure \ref{fig:CSB}).
This allows practitioners to mitigate interference within clusters while preserving statistical power and  operational feasibility.
For example, the experiment length may not be long enough relative to the mixing time to rely purely on a switchback experiment design, and the spatial network may be too densely connected to purely rely on a cluster randomized design. The application domain may also constrain the granularity of clusters we are allowed to use. By combining both spatial clusters and temporal blocks, we increase the units of randomization, and thus increase our ability to do causal estimation.

Despite the simple intuition and practical use cases, the theoretical foundations of this approach remain {\bf underexplored}.
On the surface, handling spatio-temporal interference seems straightforward, considering that (i) time can be regarded as an additional ``dimension'', and (ii) these two types of interference have been well explored separately.
However, in most work on network interference, the potential outcomes conditioned on the treatment assignments are assumed to be deterministic conditional on the treatment assignments (e.g., \citet{ugander2013graph,leung2022rate}).
In a Markovian setting, this assumption {\bf breaks down}, since past outcomes are correlated to future outcomes due to state evolution, even conditioned on the treatments.
This requires new analysis to bring together both types of interference.
The larger design space also leads to questions of optimal granularity for spatial clusters and time blocks. How much can be gained by jointly exploiting the spatial and temporal interference structure?

We consider experimentation with spatio-temporal interference on a model that encapsulates both (i) the network interference between  individuals  by a given interference graph, and (ii) the temporal interference that arises from Markovian state evolutions.
We assume that the outcome and state evolution of each individual depends solely on the treatments of their immediate neighborhood (including themselves), and that the state evolutions are independent across individuals conditioned on the treatments.

\subsection{Our Contributions }

Our main theorem shows that for a spatio-temporally clustered design, a truncated HT estimator achieves
\begin{align*}
    {\rm MSE}\lb(\wh\Delta^r\rb) =\widetilde{O}\lb(\frac{t_{\rm mix}}{T} \cdot \frac{|\{(i,i'): i\not \ind i'\}|}{2^{-2d_\Pi^{\rm max}} N^2}\rb)
\end{align*}
where $t_{\rm mix}$ is related to the mixing time of the Markovian transitions, $d_\Pi^{\rm max}$ is the maximum number of clusters neighboring an individual, and $i\not \ind i'$ if $i$ and $i'$ have neighbors within the same cluster. The first term $t_{\rm mix}/T$ matches the MSE that one expects from a pure switchback design. The second term $|\{(i,i'): i\not \ind i'\}|/2^{-2d_\Pi^{\rm max}} N^2$ is similar to the MSE that one expects from a pure graph cluster randomized (GCR) design. As a result, by combining spatial with temporal randomization, we are able to achieve a much lower MSE that benefits from both fast mixing and good clustering, e.g. growth restricted graphs \citep{ugander2013graph} or spatially derived graphs \citep{leung2022rate}.
In Table \ref{tab:our_results}, we show corollaries of our main result alongside results from the literature.
Our result bridges the literature on experimentation with spatial/network interference and temporal interference by extending the following results:
\begin{enumerate}
\item
{\bf Pure switchback experiments.} \citet{hu2022switchback} independently obtained a $\tilde O(t_{\rm mix}/T)$ MSE rate for $N=1$ using difference-in-means estimators.
Our \cref{thm:MSE} generalizes this to the $N$-individual setting using Horvitz-Thompson style estimators.
\item {\bf Network Interference.} Assuming that the interference graph satisfies the {\em $\kappa$-restricted growth} condition, \citet{ugander2023randomized}  showed that the HT estimator achieves an MSE of $\tilde O(D^2 \kappa^4/ N)$ for $T=1$ with a randomized graph cluster randomized design, where $D$ is the maximum degree.
Our \cref{coro:more_gen_than_ugander} achieves $\tilde O( t_{\rm mix}  D^3 \kappa^5 / TN)$ in the presence of Markovian temporal interference.
\end{enumerate}

Our results highlight the value of jointly randomizing over the network and time, and provides a unifying framework for selecting the optimal spatial-temporal randomized design.
In particular, when the mixing time is fast and there exists a good clustering we achieve the product of the MSE bounds from only randomizing over the network or over time.
If either the mixing time is slow or there does not exist a well separated clustering, then we can adjust the granularity of the clusters and timeblocks in our design accordingly to achieve the optimal performance between spatial and temporal randomization.

\begin{table}
\centering
\small
\begin{tabular}{p{3cm}|p{3.3cm}|p{2cm}|p{5cm}}
Interference Graph & (Spatio-) Clustering & Timeblocks & MSE \\ \midrule
No edges \newline (no interference) & one node per cluster & $r = \tilde{O}(t_{\rm mix})$ & \cref{coro:no_interference}:~$\widetilde{O}\lb(\frac{t_{\rm mix}}{ TN}\rb)$ \\ \hline
Arbitrary & all nodes in single cluster \newline (pure switchback) &  $r = \tilde{O}(t_{\rm mix})$ & \cref{coro:pure_sb}:~$\widetilde{O}\lb(\frac{t_{\rm mix}}{T}\rb)$ \newline (also \citet{hu2022switchback}) \\ \hline
Max Degree $D$ & max cluster degree $d_\Pi^{\rm max}$ \newline max cluster size $N_c$ & $r = \tilde{O}(t_{\rm mix})$ & \cref{coro:bdd_deg}:~$\widetilde{O}\lb(\frac{t_{\rm mix} D  d_\Pi^{\rm max} N_c 2^{2d_\Pi^{\rm max}}}{TN}\rb)$ \\ \hline
2D-grid spatial graph \newline edges w/in distance $h$ & uniform 2D grid &  $r = \tilde{O}(t_{\rm mix})$ & \cref{coro:spatial}:~$\widetilde{O}\lb(\frac{t_{\rm mix}h^2}{TN}\rb)$ \\  \hline
2D-grid spatial graph \newline edges w/in distance & uniform 2D grid &  $r = T$ \newline (pure GCR) & \citet{leung2022rate}:~$\widetilde{O}\lb(\frac{h^2}{N}\rb)$ \\  \hline
$\kappa$-Restricted Growth \newline Max Degree $D$ & $1$-hop-max \newline randomized clustering & $r = \tilde{O}(t_{\rm mix})$ & \cref{coro:more_gen_than_ugander}:~$\widetilde{O}\lb(\frac{t_{\rm mix} D^3\kappa^5}{T N}\rb)$ \\  \hline
$\kappa$-Restricted Growth \newline Max Degree $D$ & $1$-hop-max \newline randomized clustering & $r = T$ \newline (pure GCR) & \citet{ugander2023randomized}:~$\widetilde{O}\lb(\frac{D^2\kappa^4}{N}\rb)$ \\
\end{tabular}
\caption{{\bf Implications of Our Main Result alongside Known Results:} Our main theorem implies the following MSE bounds in several fundamental special cases.
We omit polylogarithmic terms in $N,T$ in these MSE bounds.
We also show our results alongside several known results for ``pure'' switchback and ``pure'' A/B testing under interference.
}
\label{tab:our_results}
\end{table}

\subsection{Related Work}

\noindent{\bf Violation of SUTVA.}
Experimentation is a broadly deployed learning tool in e-commerce that is simple to execute \citep{kohavi2017surprising,thomke2020building,larsen2023statistical}.
As a key challenge, the violation of the so-called {\em Stable Unit Treatment Value Assumption} (SUTVA) has been viewed as problematic in online platforms \citep{blake2014marketplace}.

Many existing works that tackle this problem assume that interference is summarized by a low-dimensional exposure mapping and that individuals are individually randomized into treatment or control by Bernoulli randomization \citep{manski2013identification,toulis2013estimation,aronow2017estimating,basse2019randomization,forastiere2021identification}.
Some work departed from unit-level randomization and introduced cluster dependence in unit-level assignments in order to improve estimator precision, including
\citealt{ugander2013graph,jagadeesan2020designs,leung2022rate,leung2023network}, just to name a few.

There is another line of work that considers the temporal interference
(or {\em carryover} effect).
Some works consider a fixed bound on the persistence of temporal interference (e.g., \citet{bojinov2023design}), while other works considered temporal interference arising from the Markovian evolution of states
\citep{glynn2020adaptive,farias2022markovian,johari2022experimental,shi2023dynamic,hu2022switchback,li2022network,li2023experimenting}.
Apart from being limited to the single-individual setting, many of these works differ from ours either by (i) focusing on alternative objectives, such as stationary outcome \citet{glynn2020adaptive}, or (ii) imposing additional assumptions, like observability of the states \citet{farias2022markovian}.

\medskip
\noindent{\bf Spatio-temporal Interference.}
Although extensively studied separately and recognized for its practical significance, experimentation under spatio-temporal interference has received relatively limited attention  previously.
Recently, \citet{ni2023design} attempted to address this problem, but their carryover effect is confined to one period.
 Another closely related work is \citet{li2022network}.
Similar to our work, they specified the spatial interference using an interference graph and modeled temporal interference by assigning an MDP to each individual.
In our model, the transition probability depends on the treatments of all neighbors (in the interference graph).
In contrast, their evolution  depends on the {\bf sum} of the outcome of direct neighbors.
Moreover, our work focuses on ATE estimation for a fixed, unknown environment, whereas they focus on the large sample asymptotics and mean-field properties.
\citet{wang2021causal} studied direct treatment effects for panel data under spatiotemporal interference, but focused on {\em asymptotic} properties instead of finite-sample bounds.
Recently, \citet{missault2025robust} considered a multi-unit switchback problem similar to ours, but they employ a ``$m$-carryover'' effect model where the potential outcome at time $t$ depends on time periods $t-m,...,t-1$ (see their Assumption 4), where $m$ is a given constant.

\begin{figure}
    \centering    \includegraphics[width=0.8\linewidth]{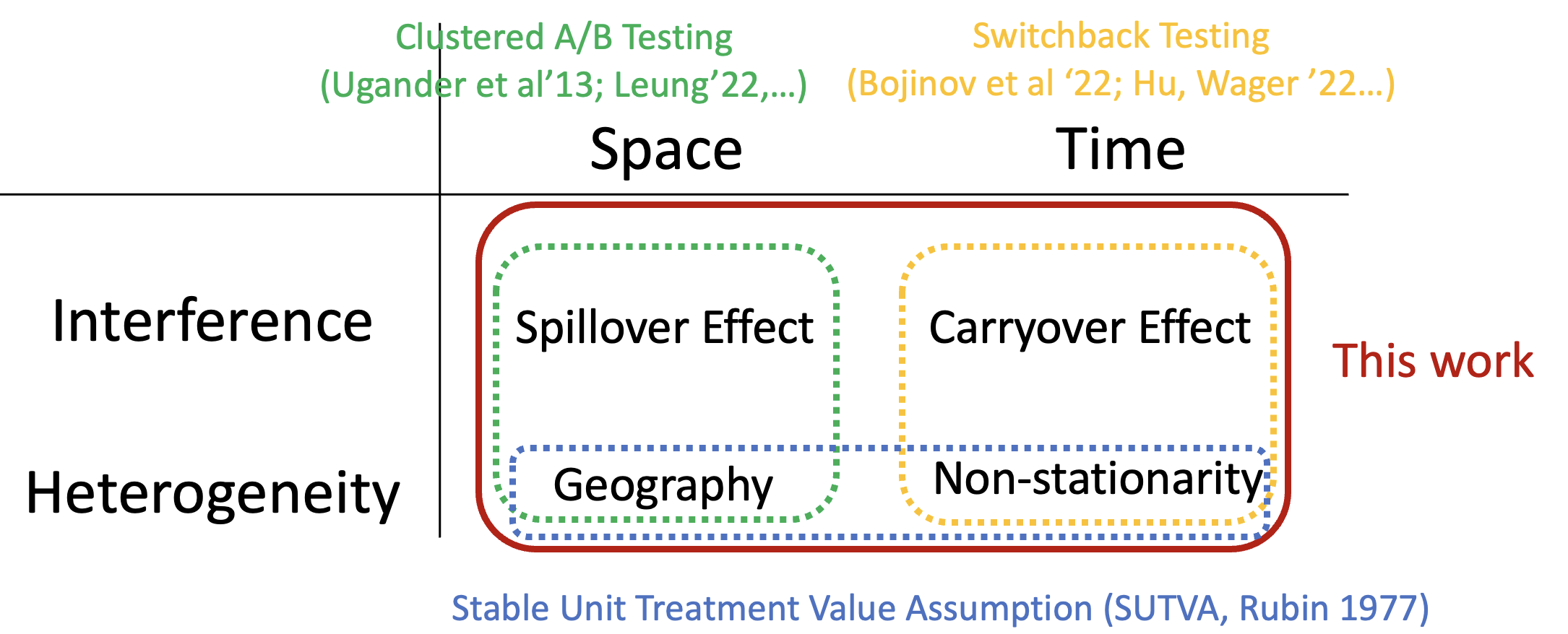}
    \caption{{\bf Positioning of this work.} By and large, a model for experimentation involves a subset of four key features, as illustrated above.
    Prior work has addressed spatial and temporal interference separately, often incorporating heterogeneity across individuals as well.
    However, in real-world applications, all four features often arise simultaneously.
    This work aims to lay the theoretical foundation for cluster switchback experiments, which are increasingly used in practice to navigate such complex interference.}
    \label{fig:positioning}
\end{figure}

\medskip
\noindent{\bf Off-Policy Evaluation (OPE)}
Since we model temporal interference using an MDP, our work is naturally related to reinforcement learning (RL).
In fact, our result on ATE estimation can be rephrased as OPE
\citep{jiang2016doubly,thomas2016data} in a multi-agent MDP: Given a {\em behavioral} policy from which the data is generated, we aim to evaluate the mean reward of a {\em target} policy.
The ATE in our work is essentially the difference in the mean reward between two target policies
(all-1 and all-0 policies), and the behavioral policy is given by clustered  randomization.
However, these works usually require certain states to be observable, which is not needed in our work.
Moreover, these works usually impose certain assumptions on the non-stationarity, which we allow to be completely arbitrary.
Finally, we focus on rather general data-generating policies (beyond fixed-treatment policies) and estimands (beyond ATE), compromising the strengths of the results.
\cref{fig:positioning} summarizes how our work relates to prior work on spatial and temporal interference.

%% file: sections/model.tex
Consider a horizon with $T$ {\em rounds} and $N$ {\em individuals}, where each individual is randomly assigned to treatment (``1'') or control (``0'') at each time. Formally, the treatment assignment is denoted by the binary matrix $W\in\{0,1\}^{N\times T}$.
We focus on {\em non-adaptive} designs where $W$ is drawn at the beginning and hence is independent of {\em all} other variables, including the individuals' states and outcomes.

\subsection{Spatial Interference and Markovian Dynamics}

We model the spatial interference between individuals using an interference graph $G=(U,E)$ where $|U| = N$ and each node represents an individual.
We refer to the direct neighbors of individual $i \in U$ as the {\it interference neighborhood} of $i$, denoted by $\N(i):=\{i\}\cup \{i\in U\mid (i,j)\in E\}$.
At each time $t\in [T]:\{1,2,\dots, T\}$, the state of {\em individual} $i \in U$ is denoted as $S_{it}\in \S$.

\begin{assumption}[Outcome Function]
\label{assump:outcome_fn}
Given an interference graph $G$,
the {\em outcome} $Y_{it}$ of individual $i$ at time $t$ satisfies $Y_{it} = \mu_{it} (S_{it}, W_{\N(i),t}) + \epsilon_{it},$
where $\N(i)$ denotes the interference neighborhood and $\epsilon_{it}$ has mean zero.
We refer to $\mu_{it}:\S\times\{0,1\}^{\N(i)} \to [0,1]$ as the {\em outcome function}.
\end{assumption}

Write $S=(S_{it})$.
We assume that $\epsilon_{it}$ have zero cross-correlation and bounded variance:

\begin{assumption}[Uncorrelated Noise]
There is a constant $\sigma$ s.t. for all $i,i',t,t'$, we have
\[
    \mathbb E[\epsilon_{it}\mid S,W]=0 \quad {\rm and} \quad \mathbb E[\epsilon_{it} \cdot\epsilon_{i't'}\mid  S,W]\leq\sigma^2 \cdot \mathbbm{1}(i=i' \ {\rm and} \ t'=t).
\]
\end{assumption}

We assume that the temporal interference arises from independent Markovian state dynamics, where the transition probabilities are a function of the treatments of the interference neighborhood.

\begin{assumption}
[Independent Transitions Conditioned on Assignment]
\label{assump:markov_transitions}
The states across different individuals evolve independently according to Markovian transition kernels that are a function of the treatments of the interference neighborhood of $i$ at the current time $t$,
\[
    \mathbb{P}(\{S_{it}\}_{i\in U} ~|~ \{S_{it'}\}_{i\in U, t' < t}, W)
= \prod_{i \in U} \mathbb{P}(S_{it} ~|~ S_{i,t-1}, W_{\N(i),t}).
\]
We denote $P_{it}^w$ the transition kernel, i.e., $P_{it}^w(S_{i,t+1}\mid S_{i,t}) := \mathbb{P}(S_{i,t+1} ~|~ S_{i,t}, W_{\N(i),t} = w)$.
\end{assumption}
We allow the transition kernels $P_{it}^w$ to vary arbitrarily across different combinations of $i,t$ and $w\in \{0,1\}^{NT}$.
Temporal interference results from the fact that the distribution of states at the current time is a function of the vector of past treatments, as the state transitions produce dependence across time.

\cref{assump:outcome_fn} and \cref{assump:markov_transitions} together capture the structure of the spatial and temporal interference, implying that the mean outcome of individual $i$ at time $t$, conditioned on the treatment vector, is only a function of the treatments of $i$'s interference neighborhood up to time $t$,
\begin{align}
\label{eq:true_exposure_map}
    \ho{E}[Y_{it}\mid W] = \ho{E}[Y_{it}\mid W_{\N(i),t'\leq t}].
\end{align}

A key assumption as introduced in \citet{hu2022switchback} that allows for estimation despite temporal interference is {\em rapid mixing}:
\begin{assumption}[Rapid Mixing]\label{assu:mixing}
There exists a constant $t_{\rm mix}>0$ such that for any individual $i\in U$, time $t\in [T]$, interference neighborhood treatment vector $w$ and distributions $f,f'$ over $\S$, we have
\[
d_{\rm TV}(f P^w_{it}, f' P^w_{it}) \le e^{- 1/ t_{\rm mix}} \cdot d_{\rm TV}(f,f').
\]
\end{assumption}

While the complete model dynamics are specified by the sequence \[W, \{(S_{i1},Y_{i1})\}_{i\in U},\dots,\{(S_{iT},Y_{iT})\}_{i\in U},\]
we only observe the treatments and outcomes $W, Y$.
In particular, we do {\bf not} observe the underlying state variables, which excludes the use of estimators that directly try to learn the Markovian dynamics.
Modeling the details of the Markovian dynamics often assumes that we correctly specified the state, which in complex applications is itself a very difficult modeling decision. By not relying on observations of the state, our results are quite robust, as we only require that there exists some state representation under which the dynamics are Markovian. While estimation will not depend on the state representation, the rapid mixing parameter $t_{\rm mix}$ may reflect the complexity of the state representation.

\begin{remark}[Extension to general exposure mapping]
While \cref{assump:outcome_fn} and \cref{assump:markov_transitions} allow for the outcome function and the transition kernels to depend on the full neighborhood treatment vector, it is natural and straightforward to impose further structure by assuming that the dependence on the neighborhood treatment vectors is mediated through a given {\em exposure mapping} $\phi_i$. The exposure mapping is a function that maps the treatment assignment of $i$'s interference neighborhood to a discrete exposure level. \cref{assump:outcome_fn} and \cref{assump:markov_transitions} would then be modified in the following way:
\begin{itemize}
    \item {\bf \cref{assump:outcome_fn} (modified for exposure mapping $\phi_i$):} The mean outcome of $i$ depends on the state of $i$, treatment of $i$, and neighborhood exposure level $\phi_i(W_{\N(i)\setminus i,t})$, formally,
    \[
        Y_{it} = \mu_{it}\lb(S_{it},W_i,\phi_i(W_{\N(i)\setminus i,t})\rb) + \epsilon_{it}.
    \]
    \item {\bf \cref{assump:markov_transitions} (modified for exposure mapping $\phi_i$):} The transition kernel that governs the dynamics of individual $i$'s state is a function of the treatments of itself and its direct neighbors through the exposure mapping $\phi_i$. Formally,
    \[
        \mathbb{P}(\{S_{it}\}_{i\in U} ~|~ \{S_{it'}\}_{i\in U, t' < t}, W)
    = \prod_{i \in U} \mathbb{P}(S_{it} ~|~ S_{i,t-1}, W_i, \phi_i(W_{\N(i)\setminus i,t})).
    \]
    Note that the individual transition kernels are then indexed by exposure levels rather than the full treatment vector, which then propagates through the remaining results.
\end{itemize}
The most general exposure mapping is the identity function, and the majority of our results are presented with the identity exposure mapping for the sake of simpler notation. A common exposure mappings considered in the literature is $\phi_i(w) = \|w\|_1$, which essentially assumes that only the number of the treated neighbors matters, and not the identity of the treated neighbors.

In the simulation section we will use the $\delta$-fractional neighborhood exposure mapping from \cite{EcklesKarrerUgander17}.

\begin{definition}\label{def:FNE}
The $\delta$-fractional neighborhood exposure mapping ($\delta$-FNE) is given by
\begin{align}\label{eqn:delta_FNE}
\phi_i(W_{\N(i)\setminus i,t}) =
\begin{cases}
1 & \text{ if } \|W_{\N(i)\setminus i,t}\|_1 \geq (1-\delta) |\N(i)\setminus i|, \\
0 & \text{ if } \|W_{\N(i)\setminus i,t}\|_1 \leq \delta |\N(i)\setminus i|, \\
W_{\N(i)\setminus i,t} & \text{ otherwise}.
\end{cases}
\end{align}
\end{definition}
\cref{assump:outcome_fn} and \cref{assump:markov_transitions} modified for $\delta$-FNE would impose that the outcome functions and transition kernels for individual $i$ at time $t$ under two treatment vectors $W$ and $W'$ are the same if $W_{it} = W'_{it}$ and if the fraction of treated neighbors in $W$ and $W'$ are both at least $(1-\delta)$ or both at most $\delta$.
\end{remark}

\subsection{GATE Estimand}

Given the observations (i.e., $W,Y$), we aim to estimate the difference between the counterfactual outcomes under continuous global deployment of treatment arm $1$ and $0$, averaged over all individuals and rounds:

\begin{definition}[Global Average Treatment Effect] Let ${\cal D}$ be any initial state distribution over ${\cal S}^N$, and denote by $\ho{E}_{\cal D}[\cdot] = \ho{E}[\cdot \mid {\bf S_0}\sim {\cal D}]$.
The {\em global average treatment effect} (GATE) is \[\Delta= \Delta_{\cal D} := \frac 1{NT} \sum_{(i,t)\in U\times [T]}\Delta_{it} \quad {\rm where}\quad \Delta_{it} = \ho{E}_{\cal D}[Y_{it}\mid W = {\bf 1}]- \ho{E}_{\cal D}[Y_{it}\mid W= {\bf 0}].\]
\end{definition}
Our results still hold if ``${\bf 1}$'' and ``${\bf 0}$'' are replaced with an arbitrary pair of {\bf fixed} treatment sequences.

If the Markov chains are rapid-mixing as defined in \cref{assu:mixing}, then \cref{prop:initial_state} shows that the impact of the initial distribution $\cal D$ is a lower-order term w.r.t. the MSE.
We will therefore suppress $\cal D$.

\begin{proposition}[Initial State Doesn't Matter]
\label{prop:initial_state}
For any distributions ${\cal D,D'}$ over ${\cal S}^N$, we have
\begin{align}\label{eqn:032225}
    |\ho{E}_{\cal D}[Y_{it}] - \ho{E}_{\cal D'}[Y_{it}]| = O\lb(e^{-t/ t_{\rm mix}}\rb), \quad \forall i,t.
\end{align}
Consequently,
\[
    |\Delta_{\cal D} - \Delta_{\cal D'}| = O\lb(\frac{t_{\rm mix} \log (NT)}{NT}\rb).
\]
\end{proposition}

\cref{prop:initial_state} implies that the error caused by misspecifying the initial distribution is $\tilde O(1/NT)$, and thus it contributes a $\tilde O(1/(NT)^2)$ term to the MSE, lower order compared to our MSE bound.

\subsection{Bernoulli Clustered Switchback Design}
We focus on clustered switchback designs, which specify a distribution for sampling the treatment vector $W$ given a fixed clustering over the network.

\bdefn[Clusters] A family $\Pi$ of subsets of $U$ is a {\em clustering} (or {\em partition}) if $C\cap C'=\emptyset$ and $\cup_{C\in \Pi} C =U$ for any $C,C'\in \Pi$.
Each set $C\in \Pi$ is called a {\em cluster}.
\edefn

We independently assign treatments to the cluster-timeblock product sets uniformly.

\begin{definition}[Bernoulli Clustered Switchback Design] \label{def:block_sb}
Let $\Pi$ be a clustering for $U$.
Uniformly partition $[T]$ into {\em timeblocks} of length $\ell>0$ (except the last one).
For each block $B\subseteq [T]$ and $C\in \Pi$, draw $A_{CB} \sim{\rm Ber}(1/2)$ independently.
Set $W_{it} = A_{CB}$ for $(i,t)\in C\times B$.
\end{definition}

We assume that the clustering is given to us, and we will present our results as function of the clustering. The design is thus parameterized by the clustering $\Pi$ and the time-block length $\ell$.

\subsection{Time Truncated Horvitz-Thompson Estimator}

We consider a class of Horvitz-Thompson (HT) \cite{horvitz1952generalization} style estimators similar to that of \citet{aronow2017estimating,leung2022rate,savje2024causal,ugander2013graph,EcklesKarrerUgander17}.
The key idea is to approximate the average outcome under global treatment (or control) by restricting to observed outcomes $Y_{it}$ that match the counterfactual outcomes under global treatment (or control). These outcomes are combined via inverse probability weights, such that their expected value is close to the global average.

Due to the Markovian state transitions, the outcome of individual $i$ at time $t$ can depend on treatments of $\mathcal{N}(i)$ in all past rounds $t' < t$ through their impact on the distribution of the state $S_{i,t-1}$.
As a result, the distribution of $Y_{it}$ is only exactly equal to the distribution of the counterfactual outcome under global treatment or control when the neighborhood units $\mathcal{N}(i)$ are assigned to treatment $a$ for all time $t' < t$.
Unfortunately, this probability is exponentially small in $t$.

However, the rapid mixing property in \cref{assu:mixing} guarantees that if the neighborhood treatment vectors of $i$ in the past $r$ rounds matches the treatment of the counterfactual global treatment (or control), then the distribution of $Y_{it}$ will be exponentially close to the distribution of the counterfactual outcome under global treatment (or control), with respect to the radius $r$.
As a result, we propose the radius-$r$ time-truncated Horvitz-Thompson estimator, which computes an inverse probability weighted average of the outcomes $Y_{it}$ on $(i,t)$ for which the neighborhood treatment vectors for $i$ in the past $r$ rounds matches the counterfactual exposure level under global treatment (or control).
This corresponds to a HT estimator with a misspecified exposure mapping with respect to time.

\bdefn[Radius-$r$ Truncated Horvitz-Thompson Estimator]
For any {\em radius} $r\ge 0$, define the {\em radius-$r$ truncated} exposure as
\begin{align}
\label{eq:def_X}
    X_{ita}^r(W) := \prod_{t'= t-r}^t \mathbbm{1}
    \lb(W_{\N(i),t'} = {\bf a}\rb)
\end{align}
for any $i,t, a$ and $W\in \{0,1\}^{N\times T}$ where ${\bf a}=(a,a,a, \dots a)$.
Define the {\em exposure probability} as
\begin{align}
\label{eq:def_p}
    p_{ita}^r = \ho{P}[X_{ita}^r = 1].
\end{align}
Denote
\begin{align}  \label{eq:def_Y_Delta}
    \wh Y_{ita}^r = \frac {X_{ita}^r} {p_{ita}^r} Y_{it} \quad {\rm and}\quad  \wh \Delta_{it}^r=\wh Y_{it1}^r -\wh Y_{it0}^r
\end{align}
for $i\in U$, $t\in [T]$ and $a\in \{0,1\}$.
The {\em Radius-$r$ Truncated Horvitz-Thompson} estimator is given by
\begin{align}
    \wh \Delta^r = \frac 1{NT} \sum_{(i,t)\in U\times [T]}\wh \Delta_{it}^r.
\end{align}
\edefn

Note that $Y_{it}$ and $X_{ita}^r$ are {\em not} independent, as they both depend on the treatments in the $r$ rounds before $t$.
This estimator is biased, but its bias decays exponentially in $r$ due to rapid mixing.
A similar spatially truncated HT estimator was proposed in the spatial interference setting \citep{leung2022rate}, and utilizes the framework of truncated (and thus misspecified) exposure mappings introduced by \citet{aronow2017estimating,savje2023causal}.

\begin{remark}
The radius-$r$ truncated exposure is {\bf misspecified} in the time dimension, since the arms from $t' < t-r$ could still impact the outcome at time $t$ through the correlation of the state distributions.
The ``true'' exposure instead involves a product from $t$ all the way to $t'$, formally,
\begin{align*}
X^{\rm True}_{ita}(W) := \prod_{t'=1}^t \mathbbm{1}
\lb(W_{\N(i),t'} = {\bf a} \rb).
\end{align*}
However, the corresponding exposure probability is exponentially low in $r/\ell$, and thus by truncating the neighborhood in the time dimension, the misspecified exposure enjoys a much higher exposure probability. This leads to a natural bias-variance tradeoff in the performance of the truncated Horvitz-Thompson estimator with respect to the choice of $r$.
Moreover, it serves as a good approximation of $X^{\rm True}_{ita}$, as the rapid-mixing property implies that the correlation across long time scales is weak, and thus limits the impact that treatments from a long time ago can have on the current outcome.
\end{remark}

In our simulations, we will also consider the H\'ajek IPW estimator \cite{hajek1971comment}, which is self-normalized and thus has significantly lower variance in the finite-sample regime.

\bdefn[Radius-$r$ Truncated H\'ajek Estimator]
For any {\em radius} $r\ge 0$,
let $X_{ita}^r$, $p_{ita}^r$, and $\wh Y_{ita}^r$ be defined by \eqref{eq:def_X} \eqref{eq:def_p}, and \eqref{eq:def_Y_Delta}.
The {\em Radius-$r$ Truncated H\'ajek} estimator is given by
\begin{align}
\widetilde \Delta^r = \frac{\sum_{(i,t)} \frac{ X_{it1}^r Y_{it}}{p_{it1}^r}}{\sum_{(i,t)} \frac{ X_{it1}^r }{p_{it1}^r}} - \frac{\sum_{(i,t)} \frac{ X_{it0}^r Y_{it}}{p_{it0}^r}}{\sum_{(i,t)} \frac{ X_{it0}^r }{p_{it0}^r}}.
\end{align}
\edefn

{In the large-sample regime, the denominator converges to $NT$, and thus the H\'ajek estimator performs similarly to the HT estimator, and is also asymptotically unbiased.
In the finite sample regime the H\'ajek estimator is biased, but has a significantly lower variance than the HT estimator when the exposure probabilities are small.
As the H\'ajek estimator is less amenable to theoretical analysis, we focus our analysis on the HT estimator, but we will show comparisons with the H\'ajek estimator in simulations.
}

\begin{remark}
If \cref{assump:outcome_fn} and \cref{assump:markov_transitions} are modified according to given spatial exposure mappings $\phi_i$, then we would define the corresponding Horvitz-Thompson and H\'ajek estimators with the respect to the modified exposure variables
\begin{align}
\label{eq:def_X_phi}
    X_{ita}^r(W) := \prod_{t'= t-r}^t \mathbbm{1}
    \lb(W_{it'} = a\rb) \prod_{t'= t-r}^t \mathbbm{1}\lb( \phi_i(W_{\N(i)\setminus i,t'}) = \phi_i({\bf a})\rb).
\end{align}
We will use this variant of the estimator with $\phi_i$ given by \cref{eqn:delta_FNE} in our simulations.
\end{remark}

\subsection{Dependence Graph}
When we analyze our proposed HT estimator, we need to consider which pairs of outcomes are correlated. We define the dependence graph as a function of both the clustering and the underlying interference network to capture the correlation structure between units in the spatial network. As the treatments are jointly assigned to all the units in a cluster, and as the outcome of $i$ depends on the treatments of its neighborhood $\mathcal{N}_i$, the outcomes of $i$ and $i'$ are correlated if any of the neighbors of $i$ and $i'$ are in the same cluster.

\bdefn[Dependence Graph]
Given a partition $\Pi$ of $U$, the {\em dependence graph} is $G_\Pi:= (U,E_\Pi)$ where for any $i,i'\in U$ (possibly identical), we include an edge $(i,i')$ in $E_\Pi$ (and denote $i\not\ind i'$) if there is a cluster $C\in \Pi$ s.t. $C\cap \N(i)\neq \emptyset$ and $C \cap \N(i')\neq \emptyset.$
\edefn

\begin{figure}
\centering\includegraphics[height=4cm,width=7cm]{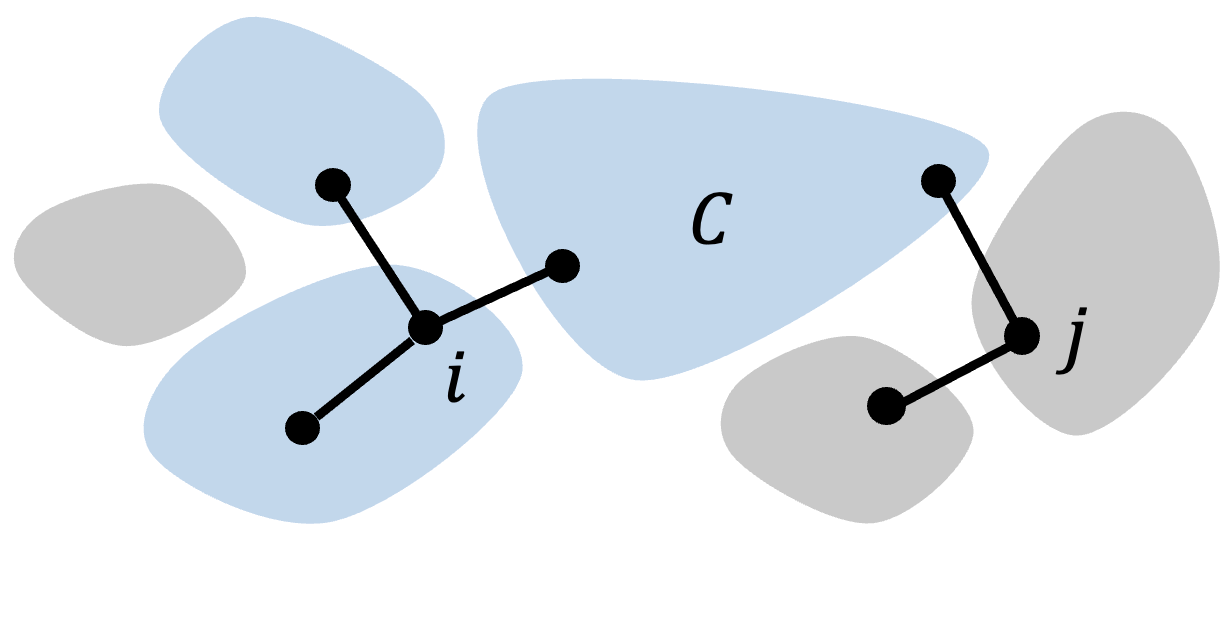}
\caption{{\bf Dependence Graph:} The regions correspond to the clusters in a partition $\Pi$.
Units $i,j$ intersect a common cluster $C$, so $(i,j)\in E_\Pi$ (or $i\not \ind j$).}
\label{fig:CIP}
\end{figure}
We illustrate this in \cref{fig:CIP} where $(i,i')$ in $E_\Pi$, denoted $i\not\ind i'$, as they both have neighbors that are members of cluster $C$. The reader should not confuse the dependence graph with the interference graph. The former is, in fact, always a supergraph of the latter.
For example, if each cluster in $\Pi$ is a singleton, then the dependence graph is the second power of the interference graph.

The dependence graph has the following useful property:
If $(i,i')\notin E_\Pi$, then $i,i'$ do not intersect any common cluster, and hence their outcomes and exposure mappings are independent under a Bernoulli clustered switchback design.

\begin{lemma}[Independence for Far-apart Individuals]
\label{prop:indep_far_apart}
Fix a partition $\Pi$ and $r\ge 0$.
Suppose $i,i'\in U$ and $(i,i')\notin E_\Pi$.
Then, for any $t,t'$, if the treatments are assigned according to a Bernoulli clustered switchback design with clustering $\Pi$, we have $\wh\Delta^r_{it} \ind \wh\Delta_{i't'}^r$.
\end{lemma}

This result is the consequence of the following facts:
(1) $\wh \Delta_{it}^r$ only depends on the treatments of the clusters that intersect $\mathcal{N}(i)$,
(2) $i\not \ind i'$ implies ${\cal C}_i \cap {\cal C}_{i'} = \emptyset$ where ${\cal C}_i$ denotes the collection of clusters that $i$ intersects, and (3) the treatments are independently assigned to each cluster.

We also define the number of clusters that an individual's interference neighborhood intersects.

\bdefn[Cluster Degree]
Given a clustering $\Pi$ of $U$, for each $i\in U$ we define \[d_\Pi(i) := |\{C\in \Pi: C\cap \N(i)\neq \emptyset\}|.\]
\edefn

%% file: sections/main_results.tex
We analyze the performance of the truncated Horvitz-Thompson estimator under the Bernoulli clustered switchback design.
Recall that $i\not\ind i'$ if $(i,i')\in E_\Pi$ where $E_\Pi$ is the set of edges of the dependence graph induced by $\Pi$, and that $p_{ita}^r$ denotes the exposure probability
\[
    p_{ita}^r = \ho{P}\lb(\bigcap_{t'= t-r}^t \lb\{W_{\N(i),t'}  = {\bf a}\rb\} \rb).
\]

\begin{theorem}[MSE Upper Bound]
\label{thm:MSE}
Fix a clustering $\Pi$ of $U$, $r\ge t_{\rm mix}$, and $\ell\ge 0$. Denote $p^{\rm min}_i = \min_{t,a} \{p_{ita}^r\}.$
Under the Bernoulli clustered switchback randomization induced by $\Pi,\ell$, we have
\begin{align}
\label{eqn:251129}
    {\rm MSE}\lb(\wh\Delta^r\rb) \le
    4 e^{-2r /t_{\rm mix}} + 8(1+\sigma^2)\lb(\frac{r+ \ell} T\rb) \lb( \frac{1}{N^2} \sum_{(i,i'): i\not \ind i'} \frac{1}{ \sqrt{p_i^{\rm min} p_{i'}^{\rm min}}}\rb).
\end{align}
\end{theorem}

\cref{thm:MSE} follows directly from the bias and variance bounds we present in \cref{prop:bias} and \cref{prop:var}.
The first term is the bound on the bias squared, and the second term is a bound on the variance.
The minimum exposure probabilities $\{p^{\rm min}_i\}_{i \in U}$ depend exponentially on the cluster degree $d_{\Pi}(i)$,
\begin{align}
\label{eq:p_min}
    p^{\rm min}_i = 2^{-d_\Pi(i) (1+\lceil r/\ell\rceil)}.
\end{align}
Each pair of $i,i'$ whose neighboring clusters intersect contributes to the variance a term that scales exponentially in the cluster degree, i.e.  $2^{(d_\Pi(i)+d_\Pi(i')) (1+\lceil r/\ell\rceil)}$. As a result, we desire clusterings that with cluster degree $d_\Pi(i)$, while still maintaining that each cluster is not too large in order to control the number of pairs $i,i'$ for which $i\not \ind i'$.

{\begin{remark}
For $\ell = r$, the optimal choice for $\ell$ and $r$ that minimizes the upper bound in \cref{thm:MSE} is
\begin{align}
\label{eq:opt_r}
    r^* = \frac{t_{\rm mix}}{2} \cdot \log\lb(\frac{1}{2(1+\sigma^2)} \lb(\frac{T}{t_{\rm mix}}\rb)\lb(\frac{2^{-2d_\Pi^{\rm max}} N^2}{|\{(i,i'): i\not \ind i'\}|}\rb)\rb) = \widetilde{\Theta}(t_{\rm mix}/2),
\end{align}
where the ``$\widetilde\Theta$'' ignores logarithmic terms. Note that the expression still grows logarithmically in $T$ and $N^2$. This expression comes from equalizing the bias and variance terms, which would give
\begin{align}
    r^* = \frac{t_{\rm mix}}{2} \cdot \log\lb(\frac{1}{2(1+\sigma^2)} \lb(\frac{T}{r^* + \ell}\rb)\lb(\frac{N^2}{\sum_{(i,i'): i\not \ind i'} (p_i^{\rm min} p_{i'}^{\rm min})^{-1/2}}\rb)\rb).
\end{align}
As the minimum exposure probability scales exponentially in $\lceil r/\ell \rceil$, the optimal choice of $\ell$ should scale as $\Theta(r)$.
When $\ell \geq r$ and the cluster degrees are upper bounded by $d_{\Pi}^{\rm max}$, then we can simplify the exposure probabilities by $p_i^{\rm min}\geq 2^{-2d_\Pi^{\rm max}}$. While there is still an $r^*$ in the denominator of the logarithmic expression, that is a lower order term.
\end{remark}

By plugging the choice of $\ell = r = r^*$ into \cref{thm:MSE} we obtain the following corollary that reveals the direct dependence of the bound on the mixing time and the properties of the graph.
\begin{corollary}
\label{coro:simplified_MSE}
For the choice of $\ell = r = r^*$ as given by \cref{eq:opt_r}, the MSE is bounded by
\begin{align*}
{\rm MSE}\lb(\wh\Delta^r\rb) &\leq 8(1+\sigma^2)\lb(\frac{t_{\rm mix}}{T}\rb) \lb( \frac{|\{(i,i'): i\not \ind i'\}|}{2^{-2d_\Pi^{\rm max}} N^2}\rb) \cdot \log\lb(\frac{e}{2(1+\sigma^2)} \lb(\frac{T}{t_{\rm mix}}\rb)\lb(\frac{2^{-2d_\Pi^{\rm max}} N^2}{|\{(i,i'): i\not \ind i'\}|}\rb)\rb), \\
&=\widetilde{O}\lb(\frac{t_{\rm mix}}{T} \cdot \frac{|\{(i,i'): i\not \ind i'\}|}{2^{-2d_\Pi^{\rm max}} N^2}\rb).
\end{align*}
\end{corollary}
}

\begin{proposition}
[Bias of the HT estimator]
\label{prop:bias}
For any $r\ge 0$, we have
\[
    \lb|\ho{E}\lb[\wh \Delta^r\rb] - \Delta\rb| \le 2 e^{-r /t_{\rm mix}}.
\]
\end{proposition}

{The bias arises from the fact that the exposure $X_{ita}^r$ used for the Horvitz-Thompson estimator only considers the past $r$ neighborhood treatment vectors, whereas technically the entire history of treatments in the past can impact the current outcome. Due to the rapid mixing assumption however, the dependence decays exponentially in time.
As a result the bias decays exponentially in $r$ at a rate governed by $t_{\rm mix}$.
This suggests that we should always choose $r \gg t_{\rm mix}$ so that the bias is small.
The proof of \cref{prop:bias} is given in \cref{sec:bias_proof}.
}

\begin{proposition}
[Variance of the HT estimator]
\label{prop:var}
Fix a clustering $\Pi$ of $U$, $r \geq t_{\rm mix}$, and $\ell\ge 0$.
Denote $p^{\rm min}_i = \min_{t,a} \{p_{ita}^r\}.$
Under the Bernoulli cluster switchback randomization induced by $\Pi,\ell$, we have
\begin{align*}
{\rm Var}\lb(\wh \Delta^r\rb) & \le 8(1+\sigma^2)\lb(\frac{r+ \ell}{T} \rb) \lb( \frac{1}{N^2} \sum_{(i,i'): i\not \ind i'} \frac{1}{ \sqrt{p_i^{\rm min} p_{i'}^{\rm min}}}\rb).
\end{align*}
\end{proposition}

The full proof of \cref{prop:var} is in \cref{sec:variance_proof}, but the key pieces are to upper bound pairwise correlations using the following:
\begin{itemize}
    \item[(1)] We bound the covariance of terms that are spatially far apart by   \cref{prop:indep_far_apart}, which states that for $(i,i')\notin E_\Pi$, denoted by {$i \ind i'$}, ${\rm Cov}(\wh\Delta_{it}^r, \wh\Delta_{i't'}^r) = 0$.
    \item[(2)] We bound the covariance of terms close in time by \cref{lem:close-by} in the appendix.
    It uses Cauchy-Schwarz to show that for all $i,i'\in U$ and $t,t'\in [T]$, ${\rm Cov}(\wh\Delta_{it}^r, \wh\Delta_{i't'}^r) \le 4(1+\sigma^2)/ \sqrt{p_{it}^r\cdot p_{i't'}^r}$.
    \item[(3)] We bound the covariance of terms far apart in time using \cref{lem:cov_bound} in the appendix. It uses the rapid-mixing property to show that for all $i,i'\in U$ and for $t,t'\in [T]$ s.t. $t+\ell < t'-r$ and $\lceil t'/\ell\rceil > \lceil t/\ell\rceil$, ${\rm Cov}(\wh\Delta_{it}^r, \wh\Delta_{i't'}^r) \le 8e^{-|t'-t|/t_{\rm mix}}.$
\end{itemize}
For the choice of $r \geq t_{\rm mix}$, it follows that the sum of the covariance terms bounded in step (3) are dominated by the sum of the covariance terms bounded in step (2), such that only the bounds from step (2) persist in the final result.

{\begin{remark}
\cref{thm:MSE} is stated for arbitrary, fixed cluster randomization.
In particular, we use the property that the treatments are assigned independently across different clusters.
For any pure-spatial randomization that partitions individuals into spatio-blocks, our result hold when a temporal component is introduced via Markovian dynamics.
Certain designs, such as the conflict graph design \cite{kandiros2024conflict}, violate the required independence.
In these cases, \cref{thm:MSE} does not apply.
We provide further discussion on extensions beyond GCR designs in \cref{subsec:RGCR}.
\end{remark}
}

\subsection{Comparison with pure spatial/temporal randomization}
As shown in \cref{coro:simplified_MSE},the MSE under a clustered switchback design with an optimal choice of $r = \widetilde{\Theta}(t_{\rm mix})$ and $\ell = \widetilde{\Theta}(t_{\rm mix})$ is approximately the product of the MSE under pure switchback and pure cluster randomized designs, which shows the gain from using a clustered switchback design that jointly randomizes over the network and time. In particular, the term $t_{\rm mix}/T$ MSE bound is typical for switchback designs, and the term $2^{2d_\Pi^{\rm max}} |\{(i,i'): i\not \ind i'\}|/N^2$ is similar to the MSE bounds typical for cluster randomized designs.

The magnitude of the gain due to jointly randomizing over the network and time depends on both the mixing time and the clustering quality.
Recall that $t_{\rm mix}$ captures how quickly the effect of a past treatment is washed away, and we will always choose $r = \widetilde\Theta(t_{\rm mix})$ in order to keep the bias small.
Large $t_{\rm mix}$ can arise when the underlying Markov process has an ``absorbing'' state — e.g., in customer-retention settings where users may churn permanently, or in high-stakes clinical trials where patients may die.
When $t_{\rm mix}$ is large, i.e. on the order of $\Theta(T)$, then the term $t_{\rm mix}/T = \widetilde\Theta(1)$.
As a result the performance as given in  \cref{thm:MSE} scales similarly to the performance of a pure cluster randomized design, and we might as well use a pure cluster randomized design, letting $r = \ell = T$.

At the other extreme, when $t_{\rm mix}$ is small but the interference graph is dense, there may not exist any good clustering which controls the cluster degrees while keeping a sparse dependence graph. As a result a pure switchback design becomes more suitable, where we simply put all nodes into a single cluster. Then the dependence graph will be a complete graph, such that $|\{(i,i'): i\not \ind i'\}|/N^2 = 1$.
As an example of such a setting, consider evaluating a ranking algorithm for a large online marketplace, where the ranking policy is refreshed frequently, so the effect of a past action washes out quickly, but the interference graph is dense as items or sellers compete within the same market.

As summarized in   \cref{tab:spatiotemp_comparison}, clustered switchback designs are well-suited in the middle-ground, when treatment effects persist long enough that pure switchbacks suffer from temporal bias, but not so long that pure spatial designs become inefficient or infeasible.
By simultaneously clustering in space and time, the clustered switchback design balances these two forces: it reduces temporal carryover effects while still controlling spatial interference.
This allows the experimenter to achieve lower mean squared error than either extreme design would provide on its own.

\begin{table}[h]
\centering
\small
\renewcommand{\arraystretch}{1.5}
\begin{tabular}{c>{\centering\arraybackslash}p{0.38\textwidth}c}
\textbf{Design} & \textbf{When to Use} & \textbf{MSE bound} \\
\midrule
Pure Switchback
& High temporal mixing; {environment stable over time} & $\tilde{O}\lb(\frac{t_{\rm mix}}{T}\rb)$ \\
Pure Spatial Clustering
& Strong spatial interference; {sparse or weakly connected graph} & $O\lb(\frac{|\{(i,i'): i\not \ind i'\}|}{2^{-2d_\Pi^{\rm max}} N^2}\rb)$\\
Clustered Switchback
& Both significant spatial interference and fast temporal dynamics & $\tilde{O}\lb(\frac{t_{\rm mix}}{T} \cdot \frac{|\{(i,i'): i\not \ind i'\}|}{2^{-2d_\Pi^{\rm max}} N^2}\rb) $\\
\end{tabular}
\caption{When to use each experimental design.} \label{tab:spatiotemp_comparison}
\end{table}

As a sanity check,  \cref{coro:no_interference} states the performance when there is no interference and we use singleton clusters.
The result follows from plugging in $|\{(i,i'): i\not \ind i'\}|=N$ and $d_\Pi^{\rm max}=1$ to   \cref{coro:simplified_MSE}.

\bcoro[No Interference]
\label{coro:no_interference}
Suppose the interference graph has no edges. For the clustering $\Pi_{\rm sgtn}$, where each cluster is a singleton, and for the choice of $\ell = r = r^*$ as given by \eqref{eq:opt_r}, the MSE is bounded by
\begin{align}
    {\rm MSE}\lb(\wh\Delta^r\rb) = O \lb(\frac{t_{\rm mix}}{TN} \cdot \log \lb(\frac{TN}{t_{\rm mix}}\rb)\rb).
\end{align}
\ecoro

\cref{coro:pure_sb} states the result for a pure switchback design when all nodes are in a single cluster. It follows from  \cref{coro:simplified_MSE}, substituting $d_\Pi^{\rm max}=1$ and $|\{(i,i'): i\not \ind i'\}|=N^2$.
\bcoro[Pure Switchback]
\label{coro:pure_sb}
Consider the clustering $\Pi_{\rm whole}$ where all individuals are in one cluster.
For the choice of $\ell = r = r^*$ as given by \eqref{eq:opt_r}, we have
\[
    {\rm MSE}(\wh\Delta^r)=O\lb( \frac{t_{\rm mix} }{T} \cdot \log\lb(\frac{T}{t_{\rm mix} }\rb)\rb).
\]
\ecoro

\begin{remark}[Comparison with \citet{hu2022switchback}]
When $N=1$, our model and design coincide with  \citet{hu2022switchback}.
They focus on a class of {\em difference-in-mean} (DIM) estimators which compute the difference in average outcomes between blocks assigned to treatment vs control, ignoring data from time points that are too close to the boundary (referred to as {\em the burn-in period}). While they show that the vanilla DIM estimators are limited to an MSE of $T^{-2/3}$, our results show that the {\em truncated} HT estimator obtains the optimal MSE, matching the improved rate of $\tilde O(t_{\rm mix}/T)$ also achieved by their bias-corrected estimator, similar to our HT estimator with $r=\ell$.
Our analysis is more general as it handles cases where $\ell \neq r$.
This is a useful distinction in practice as the block length $\ell$ is sometimes externally chosen, e.g. by an online platform, government, or nature.
By equalizing the bias and variance terms and solving for $r^*$, we can show that the optimal $r$ for a given $\ell$ satisfies
\begin{align}
    r^* = \frac{t_{\rm mix}}{2} \cdot \log\lb(\frac{2^{- d_\Pi^{\rm max}}}{2(1+\sigma^2)}\lb(\frac{T}{r^* + \ell}\rb) \lb(\frac{N^2}{|\{(i,i'): i\not \ind i'\}|} \rb) \rb) - \frac{t_{\rm mix} d_\Pi^{\rm max} \log(2)}{2}  \lb\lceil \frac{r^*}{\ell}\rb\rceil.
\end{align}
If we ignore the ceiling function, rearrange to solve for the leading term of $r^*$, it follows that
\begin{align}
r^* = \widetilde{\Theta}\lb(\frac{\ell t_{\rm mix}}{\ell + t_{\rm mix} d_\Pi^{\rm max}}\rb),
\end{align}
where the tilde notation ignores logarithmic terms. When $N=1$, then $d_{\Pi}^{\rm max} = 1$, such that $r^* \approx \min(\ell,t_{\rm mix})$.
The suppressed logarithmic terms would also show that the optimal $r^*$ increases with $T$ and the clustering quality as characterized by $N^2/|\{(i,i'): i\not \ind i'\}|$, similar to \eqref{eq:opt_r}.
\end{remark}

\begin{remark}[Relationship to \citealt{leung2022rate}]
The bias bound in \cref{prop:bias} is reminiscent of the decaying interference assumption in \citealt{leung2022rate} (albeit on distributions rather than realizations).
However, their analysis does not directly apply to our Markovian setting, since their potential outcomes are deterministic.
While the rapid mixing property implies that the temporal interference decays exponentially across time, which seems to align with Assumption 3 from \citealt{leung2022rate}, they critically assume that
$Y_{it} \ind Y_{i't'} \mid W$ for all $i,i',t,t'$,
which does not hold in our setting.
In fact, the outcomes in our Markovian setting are {\bf dependent} (albeit having weak covariance) over time, conditioned on treatment assignment.
We allow  for Markovian dependence across time within a single unit's outcomes, i.e. $Y_{it} \not \ind Y_{it'}\mid W$ for a fixed $i$ and distinct $t,t'$. We do however assume independence of the Markovian transitions across units, such that $Y_{it}\ind Y_{i't'} \mid W$ for all $t,t'$ and $i\neq i'$.
\end{remark}

\subsection{Impact of Clustering Quality}

The choice of clustering affects both the structure of the dependence graph as well as the minimum exposure probabilities.
As $\Pi$ becomes more refined, there are fewer pairs $i,i'$ with $i\not \ind i'$, but for each such pair,  the exposure probability is smaller due to its exponential dependence on the cluster degree $d_\Pi(i)$. We show special cases of our result to highlight the dependence on the network properties and clustering quality.

\begin{corollary}[Bounded-degree Graphs]
\label{coro:bdd_deg}
Let $D$ be the maximum degree of $G$, let $d_\Pi^{\rm max}$ be the maximum cluster degree associated to clustering $\Pi$, and let $N_c$ denote the maximum size of a cluster. For the choice of $\ell = r = r^*$ as given by \eqref{eq:opt_r}, the MSE is bounded by
\[
    {\rm MSE}\lb(\wh\Delta^r\rb)
= \widetilde{O}\lb(\left(\frac{t_{\rm mix}}{T}\right) \left(\frac{D  d_\Pi^{\rm max} N_c 2^{2d_\Pi^{\rm max}}}{N}\right)\rb).
\]
\end{corollary}

\begin{proof}
Note that $|\{(i,i'): i\not \ind i'\}|\le N Dd_\Pi^{\rm max} N_c$. In fact, for each $i$, there are $d_\Pi^{\rm max}$ neighboring clusters, each of which can have $N_c$ nodes, each of which can have $D$ neighbors. We plug this into \cref{coro:simplified_MSE} to obtain the bound.
\end{proof}

{A tighter upper bound on the cluster degree can be obtained for graphs where
the size of the $r$-hop neighborhoods around each node do not grow by a factor larger than $\kappa$.

\begin{definition}[$\kappa$-Restricted Growth Graphs]
A graph has the $\kappa$-restricted growth property if $|B_{r+1}(u)| \leq \kappa |B_r(u)|$ for all $u\in U$ and $r \ge 1$, where $B_{r}(u) \subseteq U$ denotes the nodes within distance $r$ from $u$.
Note that for $r= 0$, $B_r(u)$ is singleton, so $|B_1(u)|/|B_0(u)|$ can be $D+1$, which may be larger than $\kappa$.
\end{definition}

\citet{ugander2013graph} showed that when the network satisfies the $\kappa$ restricted growth property, a greedy 3-net clustering algorithm can produce a clustering where the cluster degree $d_{\Pi}^{\max}$ is bounded by $\kappa^3$ and $|\{(i,i'): i\not \ind i'\}|\le N \kappa^5 (D+1)$, see Proposition 4.2 and the proof of Proposition 4.4 in \cite{ugander2013graph}.

\

As a result, we obtain the following corollary.
As a comparison, in the static setting, \cite{ugander2013graph} bounds the variance by $D\kappa^5 2^{\kappa^3} / N$ (see proof of Proposition 4.4), which nearly matches the second term in \cref{eq:MSE_rGG}.
The extra factor of 2 in the exponent comes from the fact that full exposure requires treatments to be aligned along neighboring time blocks, s.t. $p_i^{\min}$ has an exponential dependence on $(1 + \lceil r/\ell \rceil) = 2$ for $r = \ell$.

\begin{corollary}
[$\kappa$-Restricted Growth Graphs]
\label{coro:kappa_restricted_growth}
Let $G$ be a $\kappa$-restricted growth graph with maximum degree $D$.
Let the clusters be generated by a greedy 3-net clustering algorithm.
For the choice of $\ell = r = r^*$ as given by \cref{eq:opt_r}, we have
\begin{align}
\label{eq:MSE_rGG}
    {\rm MSE}\lb(\wh\Delta^r\rb) = \widetilde{O} \lb(\left(\frac{t_{\rm mix}}{T}\right)\left(\frac{D \kappa^5 2^{2\kappa^3}}{N} \right)\rb).
\end{align}
\end{corollary}
}

We also show the results for a {\em spatially derived graph} where individuals are embedded into a $\sqrt N\times \sqrt N$ lattice.  We assume that the transitions and outcomes at a node can interfere with nodes within a hop distance $h$. In other words, we include an edge $(i,j)$ in the interference graph $G$ if $d_{\rm hop}(i,j) \le h$.

We achieve a $\tilde O(h^2 /NT)$ MSE by using a natural spatial grid clustering.
For any $s>0$, we denote by $\Pi_s$ the uniform partition of the $\sqrt N \times \sqrt N$ lattice into square-shaped clusters of size $s\times s$.
Then:

\bcoro[$h$-neighborhood Interference]
\label{coro:spatial}
Consider the interference graph over a $\sqrt N \times \sqrt N$ 2D lattice, where two nodes have an edge if $d_{\rm hop}(i,j) \le h$. For $\Pi =\Pi_{2h}$ and $\ell = r = r^*$ as given by \eqref{eq:opt_r}, the MSE is bounded by
\[{\rm MSE}(\wh\Delta^r) = \widetilde{O}\lb(\frac{t_{\rm mix} h^2}{T N} \rb)\]
\ecoro

\begin{proof}
Given the grid partition on a 2D square, and as we choose the cluster width to be twice the threshold determining edge connectivity, $d_\Pi^{\rm max} = 4$ is maximized at the corners of the clusters.
Given the spatial structure, we can argue that $|\{(i,i'): i\not \ind i'\}| = \Theta(N h^2)$, as each $i$ can have cluster degree at most 4, and these 4 clusters have at most $16 h^2$ nodes together, arranged in a $4h$ by $4h$ square. The number of nodes that are within distance $h$ of any of these nodes will still scale as $h^2$. Plugging this into   \cref{coro:simplified_MSE} results in the final bound.
\end{proof}

\subsection{Extensions Beyond Graph Cluster Randomization}\label{subsec:RGCR}

Our current theorem is stated specifically for Bernoulli clustered switchback designs, where the clustering $\Pi$ is fixed in advance, and the treatments are assigned independently across cluster-timeblocks. This excludes randomized graph cluster randomization (RGCR) \cite{ugander2023randomized}, which samples a random clustering from a distribution, and the conflict graph design \cite{kandiros2024conflict}, which assigns treatments to nodes in the order of their magnitude in the leading eigenvector of the adjacency matrix of the graph.

A modification of our analysis could be implemented to obtain analogous results that combine RGCR with a temporal switchback component. Our analysis relies on independence of the treatment assignments across time-blocks, thus to combine RGCR with a temporal switchback, we would need to independently sample a treatment assignment vector for each time-block. Importantly this means that we would independently resample a random clustering $\Pi_B$ for each timeblock $B$, and then assign treatments to each cluster $C \in \Pi_B$ independently. This maintains independence of treatment assignments across timeblocks. Under independence of the treatment across time-blocks, the result in \cref{thm:MSE} still holds, as our results only use independence of the treatment assignments across clusters in order to restrict the set of pairs $(i,i')$ such that $\wh\Delta^r_{it} \ind \wh\Delta_{i't'}^r$.

The corresponding Horvitz-Thompson estimator for RGCR uses the marginal exposure probability with respect to the additional randomness in the clustering in the denominator, which is lower bounded by the probability that the node is in the interior of a cluster. \citet{ugander2023randomized} proposed the 1-hop max clustering, which involves first generating a random seed $Z_i\sim U(0,1)$ for each $i$ independently, then assigning $i$ to the cluster indexed by the individual in its interference neighborhood with the largest random seed, $\arg\max_{j\in \N(i)} Z_j$. They show that each node has at least $1/(D+1) \kappa$ probability to be in the interior of a cluster, which they use to lower bound the exposure probability (Theorem 4.2 in \citet{ugander2023randomized}). They additionally show that the outcomes and exposure levels of nodes that are more than distance 4 apart from each other in the graph are independent (Lemma B.6 in \citet{ugander2023randomized}).

\bcoro[Restricted-Growth Graphs]
\label{coro:more_gen_than_ugander}
Suppose $G$ satisfies the $\kappa$-restricted growth property and has maximum degree $D$. If we assign treatments independently across timeblocks according to the 1-hop-max RGCR design in \citet{ugander2023randomized}, then
\begin{align}\label{eq:MSE_RGCR}
    {\rm MSE}\lb(\wh\Delta^r\rb)
    \leq  4 e^{-2r /t_{\rm mix}} + 8(1+\sigma^2)\lb(\frac{r+ \ell}{T} \rb) \lb( \frac{(1+D)\kappa^3}{N} \rb) \lb(2(D+1) \kappa \rb)^{(1+ \lceil r/\ell \rceil)}.
\end{align}
For an appropriate choice of $\ell = r$, the MSE simplifies to
\begin{align} \label{eq:simp_MSE_RGCR}
    {\rm MSE}\lb(\wh\Delta^r\rb)
    =\widetilde{O}\left( \lb(\frac{t_{\rm mix}}{T} \rb) \lb( \frac{D^3\kappa^5}{N} \rb)\right).
\end{align}
\ecoro

\begin{proof}
Using Theorem 4.2 from \citet{ugander2023randomized} and the independence of treatment assignments across timeblocks, it follows that
\[
    p^{\rm min}_i \ge \lb(\frac{1}{2(D+1) \kappa} \rb)^{(1+ \lceil r/\ell \rceil)}.
\]
Using Lemma B.6 in \citet{ugander2023randomized} and the fact that for a $\kappa$-growth restricted graph the number of nodes that are within distance 4 of any $i$ is upper bounded by $(1+D)\kappa^3$, it follows that
\begin{align*}
    {\rm MSE}\lb(\wh\Delta^r\rb)
    \le 4 e^{-2r /t_{\rm mix}} + 8(1+\sigma^2)\lb(\frac{r+ \ell}{T} \rb) \lb(\frac{(1+D)\kappa^3}{N} \rb) \lb(2(D+1) \kappa \rb)^{(1+ \lceil r/\ell \rceil)}.
\end{align*}
The simplified MSE results from choosing
\[
    \ell = r = \frac{t_{\rm mix}}{2} \cdot \log\lb(\frac{1}{8(1+\sigma^2)} \lb(\frac{T}{t_{\rm mix}}\rb) \lb(\frac{N}{(1+D)^3\kappa^5}\rb)\rb).
\]
\end{proof}

{As compared to \cref{coro:kappa_restricted_growth}, by randomizing over clusters, we can obtain significantly lower variance that scales only polynomially rather than exponentially in $\kappa$, replicating the insight from \cite{ugander2023randomized}. Instead of paying the worst-case exposure penalty of a fixed clustering, the randomized 1-hop-max clustering guarantees each node is in the interior of a cluster with probability $1/(D+1)\kappa$, improving the scaling of $p_i^{\min}$ from depending exponentially on $\kappa^3$ to only depending polynomially on $\kappa$.

In the static setting, the variance bound shown in Theorem 4.7 of \cite{ugander2023randomized} scales as $D^2 \kappa^4 / N$, which nearly matches the second term in \cref{eq:simp_MSE_RGCR}. Our bound has an additional factor of $D$ and $\kappa$ due to the $\lceil r/\ell \rceil$ term in the exponent of \cref{eq:MSE_RGCR}, which arises from the fact that full exposure at the boundary of time-blocks requires that neighboring time-blocks have the same treatment. In addition our bound has an additional factor of $\frac{t_{\rm mix}}T$, which gains us significant variance reduction for large $T$. At a high level, our analysis separates out the reduction in variance due to the temporal component as opposed to the spatial network.}

\begin{remark}
Analyzing the conflict graph design for the standard Horvitz-Thompson estimator is significantly harder due to the dense correlations. \cite{kandiros2024conflict} handle this by instead proposing a modified Horvitz-Thompson estimator which allows for sparse correlations.
{We leave this formal extension for future work, but we believe that a similar exercise could be done to extend their result to incorporate temporal dynamics. Firstly, the conflict graph design would need to be resampled independently for each time-block to maintain independence across time-blocks.} Additionally we would need to use the modified estimator defined in section 3.4 of \cite{kandiros2024conflict} along with a time-truncated modified exposure mapping. The fast mixing property would be used to show that the correlation between terms decays exponentially fast in time, with the rate governed by the mixing time. For terms that are close in time, we would use the lemmas from \cite{kandiros2024conflict} to show that the correlations are spatially sparse. Combining these pieces together would result in a bound that is the product of a term relating to the temporal dynamics and a term relating to the spatial interference.
\end{remark}

\subsection{Discussion of exposure probabilities}
\label{sec:disc_expp}

Depending on the interference network structure, there may not exist a clustering that separates the network well into many loosely connected clusters. When a high quality clustering is difficult to find, then the exponential dependence on $d_\Pi(i)$ in \cref{eq:p_min} can still result in exposure probabilities that are too low to be practical. For example,  if $\ell=r$ and $d_\Pi(i)=3$, then $p_i^{\rm min} = 2^{-6} \approx 0.0156$, such that we will only use a very small fraction of the data. If we use a more restrictive exposure mapping that imposes structure on the role of the interference neighborhood, then the corresponding exposure probability may be significantly larger.

\begin{remark}[Smoothed Exposure Mapping - $\delta$-fractional neighborhood exposure]

If \cref{assump:outcome_fn} and \cref{assump:markov_transitions} are modified according to given spatial exposure mappings $\phi_i$, \cref{thm:MSE} still applies, where the only difference is that the exposure probability would instead be
\[
    p_{ita}^r = \ho{P}\lb(\bigcap_{t'= t-r}^t \lb\{W_{i,t'}  = a\rb\} \bigcap_{t'= t-r}^t \lb\{\phi_i\lb(W_{\N(i)\setminus i,t'}\rb)  = \phi_i({\bf a})\rb\} \rb).
\]
\end{remark}

For example, consider the $\delta$-fractional neighborhood exposure mapping in \cref{def:FNE}, with the accordingly adjusted Horvitz-Thompson estimator where the exposure variables are defined as in \eqref{eq:def_X_phi}. For $\ell=r$ and $\delta = 0.25$, then $X_{it1}^r(W) = 1$ if $W_i = 1$ and at least $3/4$ of the other neighboring nodes of $i$ are assigned also to treatment. Suppose we have a node $i$ with 8 other neighbors not including itself, 4 of which are in the same cluster as node $i$, and the other 4 nodes are evenly split between two other neighboring clusters. The cluster degree $d_\Pi(i)= 3$ and $|\mathcal{N}(i) \setminus i| = 8$. In order to satisfy the condition that at least $3/4$ of the neighbors are assigned to treatment, we can only allow for one of the other clusters to be assigned to control. As a result, $p_i^{\rm min} = (3 \cdot 2^{-3})^2 \approx 0.1406$, which is nine times larger than the minimum exposure probability without imposing structure on the exposure mapping.

%% file: sections/xp.tex
We illustrate our estimator and design via a case study, similar to the rental marketplace scenario in \cite{farias2022markovian}.
Consider a hotel booking platform such as TripAdvisor.
Suppose the platform is experimenting with increasing the resolution of images uploaded to its website.
The increased resolution may increase customer interest and booking rates, in exchange for higher datacenter costs.

To determine whether the causal effect would be worth the extra cost, the platform would like to conduct a randomized control trial.
A hotel assigned to treatment would be allowed to upload high resolution images, whereas hotels assigned to control would be limited as usual to low resolution images.
In addition to a direct effect of increasing booking rates on the treated hotel, we expect that due to competition, there would be an indirect effect of {\bf decreasing} booking rates on neighboring hotels.

\subsection{Hotel Occupancy Problem}

We model the occupancy of each hotel $i$ as a random walk where the occupancy of the hotel at time $t$ is modeled as a state $S_{it}$ that evolves as a Markov chain.
The occupancy is upper bounded by the hotel's capacity $m_i$.
With probability $\lambda^l_{it} \in [0,1]$, the state does not change (no arrivals or departures).
The 1D-Markov chain is depicted below:
\begin{center}
\begin{tikzpicture}[
    >=Stealth,
    node distance=2cm,
    state/.style={
        circle,
        draw=black,
        thick,
        minimum size=9mm,
        inner sep=0pt
    }
]

\node[state] (0) {$0$};
\node[state, right=of 0] (1) {$1$};
\node[state, right=of 1] (2) {$2$};
\node[right=1.4cm of 2] (dots) {$\cdots$};
\node[state, right=1.4cm of dots] (m-1) {$m_i\!-\!1$};
\node[state, right=of m-1] (m) {$m_i$};

\draw[->, thick, bend left=20] (0) to node[above] {$p^a_{it}$} (1);
\draw[->, thick, bend left=20] (1) to node[above] {$p^a_{it}$} (2);
\draw[->, thick, bend left=20] (m-1) to node[above] {$p^a_{it}$} (m);

\draw[->, thick, bend left=20] (1) to node[below] {$p^d_{it}$} (0);
\draw[->, thick, bend left=20] (2) to node[below] {$p^d_{it}$} (1);
\draw[->, thick, bend left=20] (m) to node[below] {$p^d_{it}$} (m-1);

\draw[->, thick] (0) to[loop above] node {$1 - p^a_{it}$} (0);
\draw[->, thick] (1) to[loop above] node {$\lambda^l_{it}$} (1);
\draw[->, thick] (2) to[loop above] node {$\lambda^l_{it}$} (2);
\draw[->, thick] (m-1) to[loop above] node {$\lambda^l_{it}$} (m-1);
\draw[->, thick] (m) to[loop above] node {$1 - p^d_{it}$} (m);
\end{tikzpicture}
\end{center}
We next explain the details.

\noindent{\bf Transition Probability: Interior States.}
Conditioned on that the state changes, the departure and arrival rates are governed by parameters $\lambda^d_{it} \in [0,1]$ and
\[
    \lambda^a_{it} = \sigma (\alpha_{it} + \beta_{it} W_{it} - \gamma_{it} C_{it})
    \quad {\rm where} \quad
     C_{it} = \min\left\{0.25, \ \max\left(0.75, \frac{\sum_{j \in \mathcal{N}_i \setminus {i}} W_{jt}}{|\mathcal{N}_i \setminus {i}|}\right)\right\}
\]
with $\sigma(\cdot)$ denoting the sigmoid function.
For $s \in \{1, 2, \dots m_i - 1\}$, the transition probabilities are
\begin{align*}
    \mathbb{P}(S_{i,t+1} = s+1 ~|~ S_{it} = s) &= \frac{(1 - \lambda^l_{it})\lambda^a_{it}}{\lambda^d_{it} + \lambda^a_{it}} =: p^a_{it}, \\
    \mathbb{P}(S_{i,t+1} = s-1 ~|~ S_{it} = s) &= \frac{(1 - \lambda^l_{it})\lambda^d_{it}}{\lambda^d_{it} + \lambda^a_{it}} =: p^d_{it}, \\
    \mathbb{P}(S_{i,t+1} = s ~|~ S_{it} = s) &= 1 - p^a_{it} - p^d_{it} =: \lambda^l_{it},
\end{align*}
and
\[
    \mathbb{P}(S_{i,t+1} = s' |\ S_{it} = s)=0
\]
for any $s' \notin \{s,s-1,s+1\}$,
where $p^a_{it}$ and $p^d_{it}$ denote the probability of an arrival or a departure.

\noindent{\bf Transition Probability: Boundary States.}
For states $s = 0$ and $s = m_i$ at the boundary, we have
\begin{align*}
    \mathbb{P}(S_{i,t+1} = 1 ~|~ S_{it} = 0) = p^a_{it},
    \quad&\mathbb{P}(S_{i,t+1} = 0 ~|~ S_{it} = 0) = 1-p^a_{it}, \\
    \mathbb{P}(S_{i,t+1} = m_i-1 ~|~ S_{it} = m_i) = p^d_{it},
    \quad&\mathbb{P}(S_{i,t+1} = m_i ~|~ S_{it} = m_i) = 1-p^d_{it}.
\end{align*}

\noindent{\bf Outcomes.}
The outcome $Y_{it}$, representing the platform's profit, involves a fixed fee that the hotel pays to be listed on the platform, and a per-occupancy commission fee for each occupant.
More precisely,
\[
    Y_{it} = a_{i} + b_{it} S_{it}.
\]
The hotel capacities $m_i$ are known but the transition and reward parameters $\lambda^l_{it}, \lambda^d_{it}, \alpha_{it}, \beta_{it}, \gamma_{it}, a_{i}, b_{it}$ are unknown.
We generate the parameters by sampling each independently as
\[
    \lambda^l_{it} = 0.1,\ \lambda_{it}^d \sim U[0,1], \ \alpha_{it} \sim N(-5,2), \ \beta_{it} \sim N(10,4), \gamma_{it} \sim U[0,0.5 \beta_{it}], \ a_{i} \sim U[0,1],\ b_{it} \sim U[1, 1.2]
\]
For each hotel $i$, the occupancy at time 0 is sampled uniformly from $[0,m_i]$.

\begin{figure}[thbp]
    \centering    \includegraphics[width=0.6\textwidth]{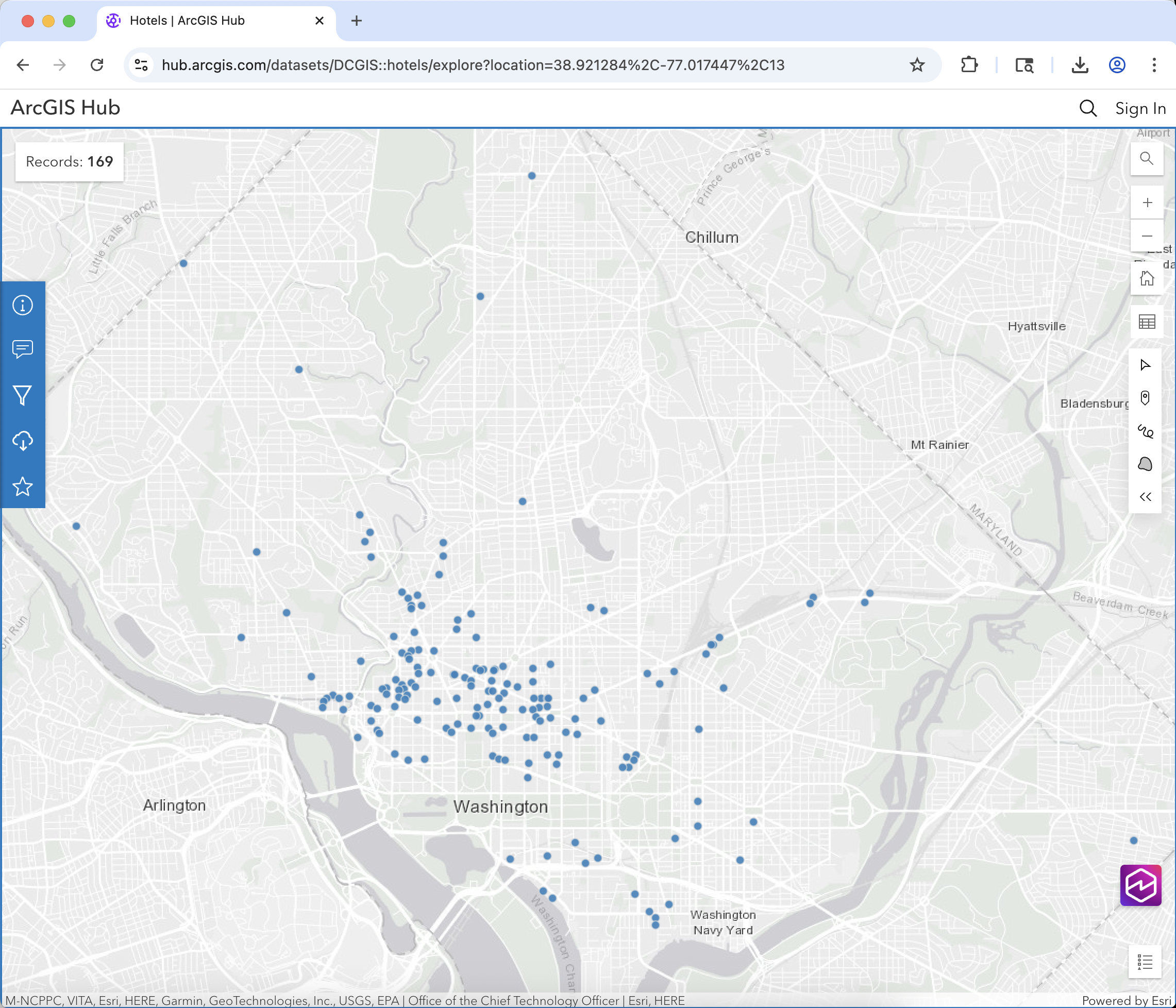}
    \caption{Locations of 169 Hotels in Washington DC provided in dataset \cite{dc_hotels_2025}.}
    \label{fig:hotel_map}
\end{figure}

\noindent{\bf Interference.}
The network interference arises from the fact that the transitions depend on the competition $C_{it}$ which represents the fraction of treated neighbors of $i$ (excluding itself) at time $t$.
Since the competition term takes the fraction of treated neighbors and projects it to be within 0.25 and 0.75, there is only network interference if at least 0.25 fraction of the neighbors are treated, and if at least 0.75 of the neighbors are treated the dynamics are equivalent to if all the neighbors are treated. This reflects the intuition that if a large majority of the neighbors are treated or control then the effective interference is the same as under global treatment or global control.

While the transition probabilities across hotels are correlated through the treatment status of neighboring hotels, we assume that the states at different hotels evolve {\bf independently} of each other.
We leave the extension to correlated transitions as an interesting future direction.
This requires modeling the joint distribution over individual hotel state transitions as a single joint Markov chain over the entire vector of hotel occupancies, significantly increasing the complexity of the correlation structure.

\subsection{Networks}
We consider two settings.
In the first, we use a synthetic network in which the locations of $n=100$ hotels are sampled uniformly over the unit square.
Two hotels are neighbors if their Euclidean distance is less than $\kappa=0.2$.
The hotel capacities are set to $m=30$ uniformly.

In the second, we use a real-world network constructed from a dataset of 169 hotels in Washington D.C. \cite{dc_hotels_2025} (see \cref{fig:hotel_map}).
The dataset contains locations and attributes of hotels, created as part of the D.C. Geographic Information System for the D.C. Office of the Chief Technology Officer and participating D.C. government agencies.
We used the x and y coordinates of each hotel, along with an attribute stating the capacity (maximum number of rooms) in each hotel.\footnote{We divided the reported capacities by 10 to limit the size of the state space for the purpose of computation}
The capacity was only provided for 141 hotels.
For the remaining 28 hotels, we randomly sampled with replacement artificial capacities from the 141 provided datapoints to maintain a similar empirical distribution.
We scaled the coordinates to be contained in the unit square, setting pairs of hotels to be considered neighbors of each other if their euclidean distance is less than $\kappa = 0.035$, which corresponds to 10-minute walking distance.

\begin{figure}[tb]
  \centering
  \begin{minipage}[t]{0.48\textwidth}
    \centering
    \includegraphics[width=\linewidth]{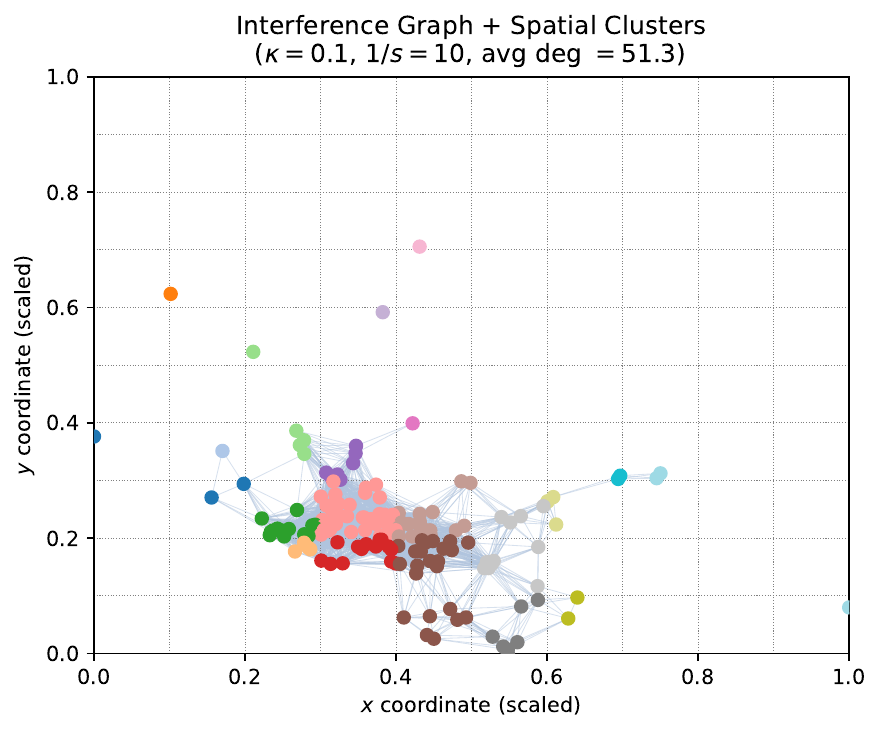}
    \caption{Interference graph of D.C. hotels ($\kappa=0.1$, $1/s=10$). Edges connect hotels within distance $\kappa$.
    Colors represent the spatial clustering.}
    \label{fig:hotel_graph}
  \end{minipage}
  \hfill
  \begin{minipage}[t]{0.48\textwidth}
    \centering
    \includegraphics[width=\linewidth]{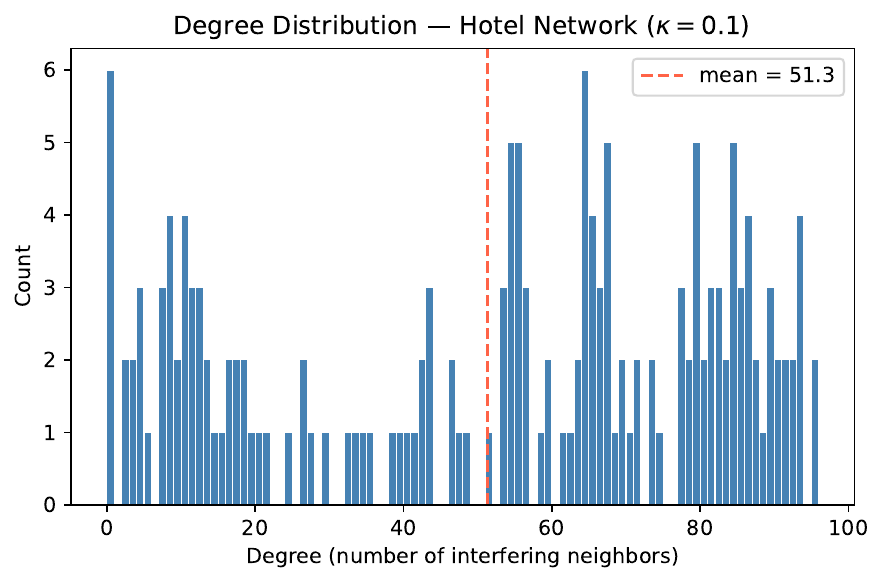}
    \caption{Degree distribution of the interference graph.
    The distribution is right-skewed, reflecting the non-uniform density.}
    \label{fig:hotel_degree}
  \end{minipage}
\end{figure}

The D.C. hotel network consists of $n = 169$ hotels whose coordinates are
scaled to the unit square $[0,1]^2$.
Two hotels are connected by an edge if their Euclidean distance is at most $\kappa = 0.1$; summary statistics are reported in
\cref{tab:network_stats}, and the interference graph and its degree distribution are shown in \cref{fig:hotel_graph,fig:hotel_degree}.
The degree distribution is highly heterogeneous, reflecting the spatial
concentration of hotels in the downtown core, with a small fraction of
peripheral hotels having no interfering neighbors.

\begin{table}[htb]
\centering
\small
\begin{tabular}{lr}
Statistic & Value \\
\midrule
Nodes & 169   \\
Edges                           & 4{,}332 \\
Mean degree                     & 51.3  \\
Median degree                   & 56    \\
Std.\ degree                    & 30.3  \\
Min / Max degree                & 0 / 95 \\
10th / 90th pct.\ degree        & 8 / 87 \\
Isolated nodes                  & 6 (3.6\%) \\
Connected components            & 7     \\
\end{tabular}
\caption{Summary statistics of the DC hotel interference graph
($\kappa=0.1$).}
\label{tab:network_stats}
\end{table}

We use a cluster switchback design which has a hyperparameter $\ell$ determining the time-block length, and a hyperparameter $s$ which defines the width of the spatial clusters.
By setting $\ell = T$ or $s = 1$ we recover a pure switchback or a pure spatial cluster randomized design.

We use a Horvitz-Thompson estimator with misspecified exposure mapping, where $r$ is the radius/recency used to truncate the exposure mapping, and $\delta$ is used to relax the exposure mapping.
In our example, for $\delta \leq 0.25$ the spatial exposure mapping is correctly specified.
We will also consider the H\'ajek estimator \cite{hajek1971comment}, which is self-normalized  version of HT estimator, and thus is expected to have significantly lower variance in the finite sample regime.

\subsection{Simulation Results: Synthetic Network}
\label{sec:simulation}

\begin{figure}[ht]
  \centering
  \includegraphics[width=0.7\linewidth]{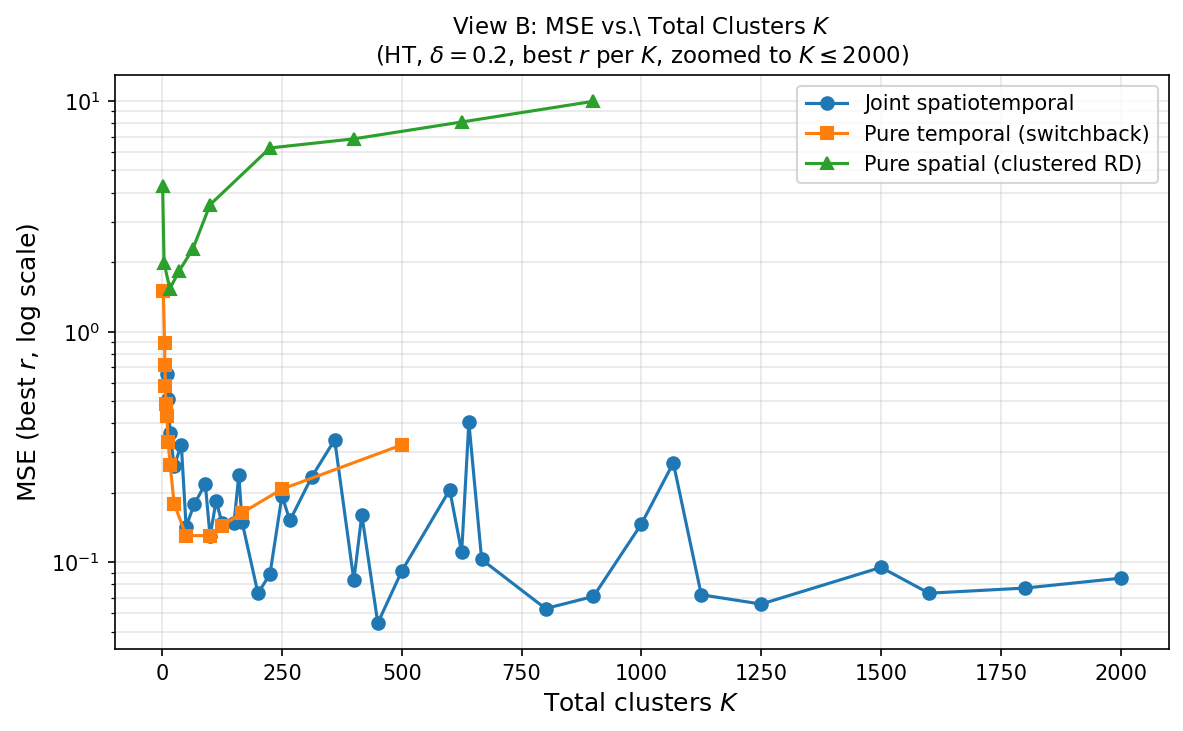}
  \caption{MSE vs. total (spatio-temporal) clusters $K$, with HT estimator, $\delta=0.2$. Joint spatio-temporal design has lower MSE than both pure designs.
  }
  \label{fig:design_viewB}
\end{figure}

We evaluate the finite-sample performance of the proposed estimators on a synthetic network with $n = 100$ spatial units distributed uniformly on $[0,1]^2$, $T = 10{,}000$ time periods.
Throughout we fix $\delta = 0.2$ for the fractional neighborhood condition.
We compare three designs - joint spatiotemporal, pure temporal (switchback), and pure spatial (clustered RD) - and three estimators: Horvitz-Thompson (HT),  Haj\'{e}k, and the
difference-in-means (DM) estimator with and without a burn-in period.

\medskip
\noindent{\bf Design Comparison.} We compare the three designs in terms of the best achievable MSE of the HT estimator, optimized over all tuning parameters, i.e., size $s$ of spatial cluster, time block length $\ell$, and recency $r$ for HT/Haj\'{e}k estimator.
\cref{fig:design_viewB} fixes the total number $K = (1/s)^2 \cdot (T/\ell)$  of spatiotemporal clusters, comparing designs at equal complexity.
The joint design remains superior, showing that the gain is not merely due to using more clusters.
Numerically, the best achievable MSE for HT estimator is $0.054$ for the joint design, $0.130$ for pure temporal ($2.39\times$ higher), and $1.534$ for pure spatial
($28.2\times$ higher).

\begin{figure}[htbp]
  \centering
  \includegraphics[width=\linewidth]{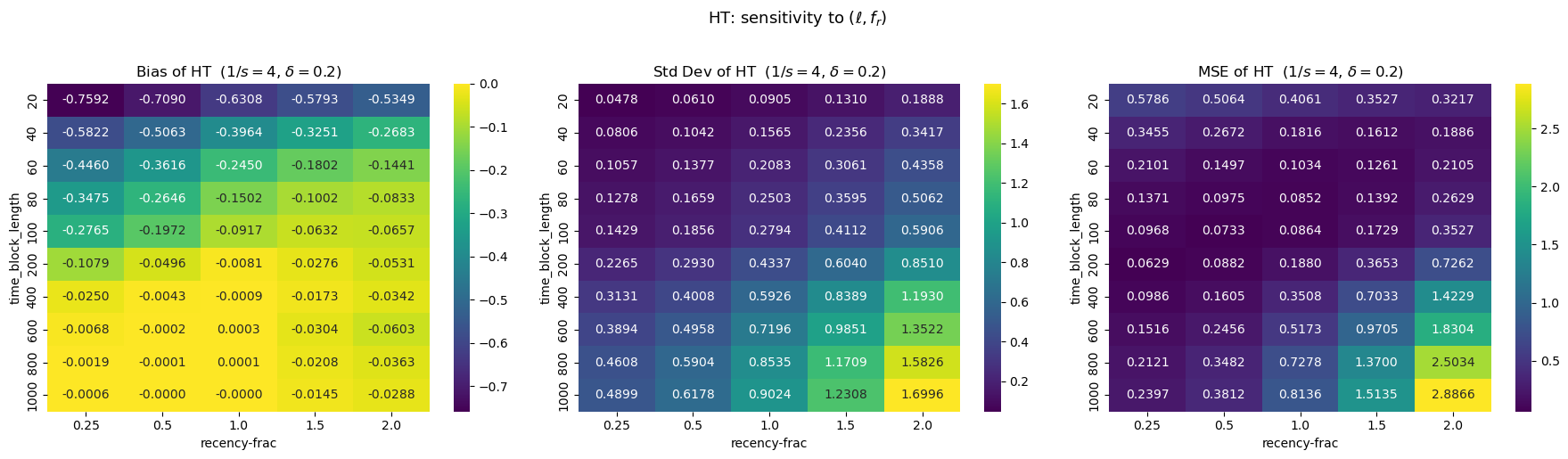}
  \caption{%
    Heatmaps of bias (left), standard deviation (centre), and MSE (right)
    for the HT estimator as a function of time-block length $\ell$ (rows)
    and fractional recency $f_r = r/\ell$ (columns), fixing $1/s = 4$ and
    $\delta = 0.2$ (Uniform network, $n=100$, $T=10{,}000$).
    MSE is U-shaped in $f_r$ for every $\ell$, with the optimum near
    $f_r^* \approx 0.5$.
  }
  \label{fig:christina_fig1}
\end{figure}

\medskip
\noindent{\bf Sensitivity to $(\ell, f_r)$.}
\cref{fig:christina_fig1} shows the bias, standard deviation, and MSE of the HT estimator as a function of the time-block length $\ell$ (rows) and the fractional recency $f_r = r/\ell$ (columns), fixing $1/s = 4$ and
$\delta = 0.2$.
Two complementary trends are noticeable.
First, for a fixed $\ell$, the MSE is U-shaped in $f_r$: small $f_r$ (short recency window) yields large bias because carry-over interference from the previous block is not adequately discarded, while large $f_r$ yields large variance because very few observations satisfy the recency criterion in the HT estimator (for not being discarded).
The optimal $f_r^*$ lies in the range $[0.5, 1.0]$ across all $\ell$ values shown, suggesting that setting $r \approx \ell/2$ is a robust practical choice for this network.
Second, for a fixed $f_r$, increasing $\ell$ initially reduces MSE (fewer blocks means less boundary ``contamination''), but eventually increases it once $T/\ell$ becomes too small to average out Monte Carlo noise.

\begin{figure}[htbp]
  \centering
  \includegraphics[width=\linewidth]{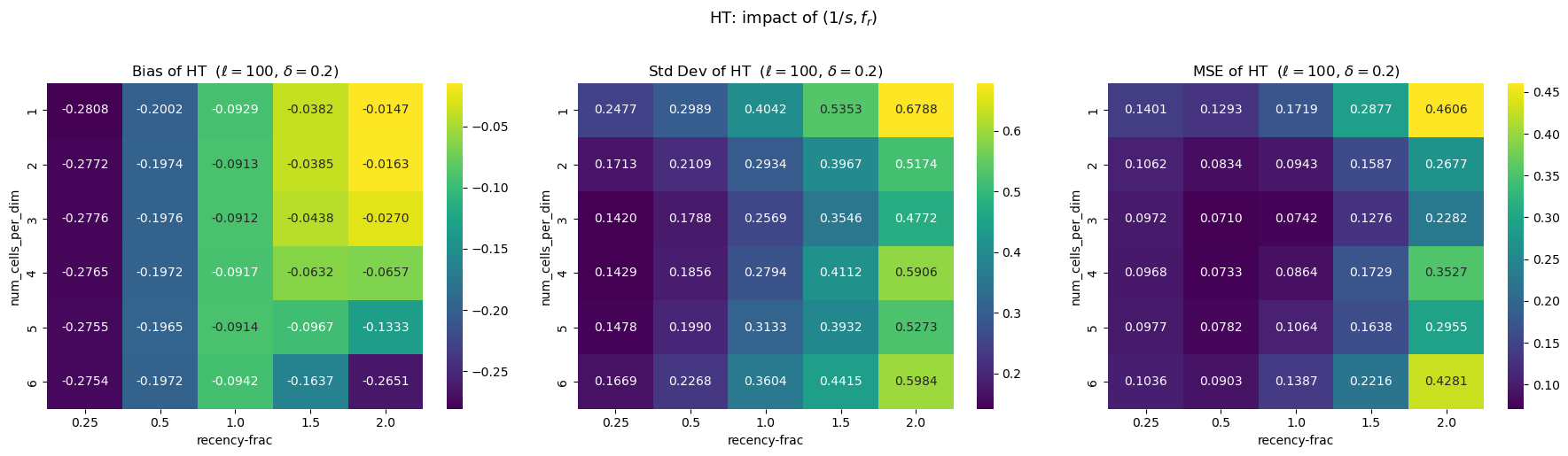}
  \caption{%
    Heatmaps of bias (left), standard deviation (center), and MSE (right)
    for the HT estimator as a function of spatial cells per dimension $1/s$
    (rows) and fractional recency $f_r$ (columns), fixing $\ell = 100$ and
    $\delta = 0.2$ (Uniform network).
    MSE is minimized at an interior value $1/s \in \{3, 4\}$; beyond this,
    variance rises as each spatial cluster contains fewer units.
  }
  \label{fig:christina_fig2}
\end{figure}

\noindent{\bf Impact of spatial granularity $s$ on HT.}
\cref{fig:christina_fig2} fixes $\ell = 100$ and shows how the HT estimator's
bias, standard deviation, and MSE vary with the number of spatial cells per dimension $1/s$ (rows) and $f_r$ (columns).
Increasing $1/s$ from $1$ to $4$ substantially reduces both bias and MSE, consistent with the theoretical prediction that finer spatial clustering reduces the spatial interference leakage across cluster boundaries.
Beyond $1/s = 4$, however, the gains diminish and variance begins to rise, because each spatial cluster contains fewer units.
This reveals an interior optimum in $1/s$: for the Uniform network with $n = 100$ and interference radius $\kappa = 0.2$, the sweet spot is around $1/s \in \{3, 4\}$.
The pattern is consistent across all $f_r$ values, indicating that spatial granularity is more consequential than the choice of recency fraction.

\begin{figure}[htb]
  \centering
  \includegraphics[width=\linewidth]{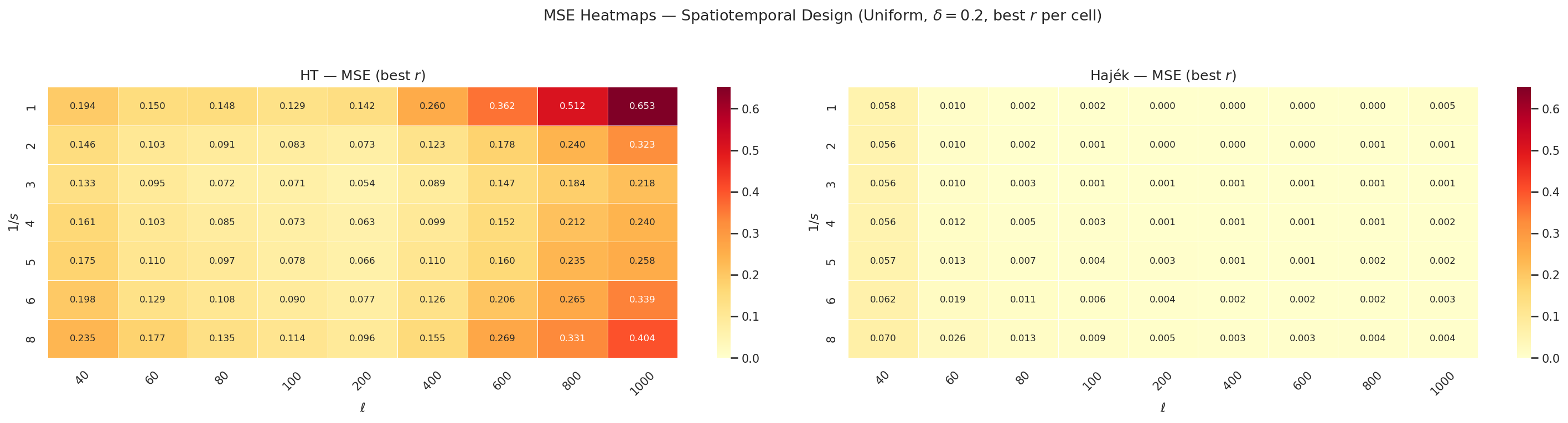}
  \caption{Heatmaps of bias, standard deviation, and MSE for the HT estimator over $(\ell, 1/s)$ ($\delta=0.2$, best $r$ per cell).
  Rows: $1/s \in \{1,\ldots,8\}$; columns: $\ell \in \{40,\ldots,1000\}$.
  This is expected, since when $\ell$ is small, ``contaminated'' data (rounds) make up a higher proportion of the entire data.}
  \label{fig:heatmap_HT}
\end{figure}

\medskip
\noindent{\bf Parameter sensitivity: heatmaps over $(\ell, s)$.}
\cref{fig:heatmap_HT} shows the bias, standard deviation, and MSE over the joint parameter space $(\ell, s)$ for the HT and Haj\'{e}k estimators, with the recency $r$ chosen to minimize MSE for each combination.
The HT heatmap shows that MSE is minimized at intermediate $\ell$ (around $200$-$400$) and moderate $1/s$ (around $3$-$5$), consistent with the bias-variance trade-off in both dimensions.
The Haj\'{e}k heatmap shows the same qualitative pattern but with
substantially lower MSE throughout.

\begin{figure}[htbp]
  \centering
  \includegraphics[width=\linewidth]{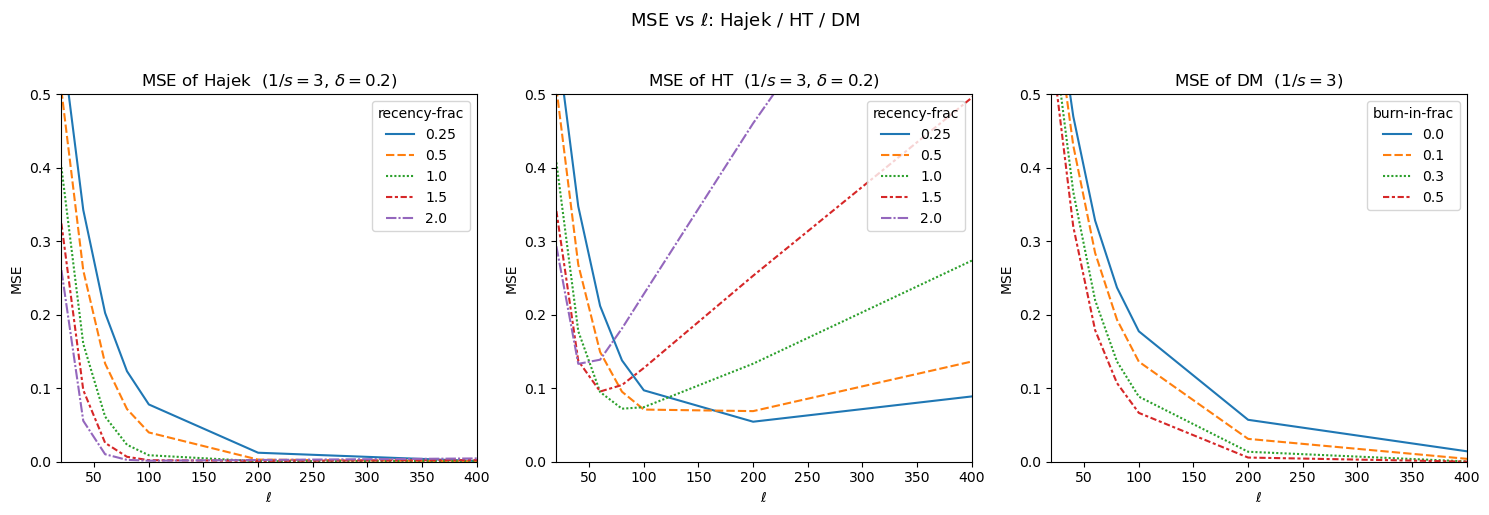}
  \caption{MSE vs.\ time-block length $\ell$ for the Haj\'{e}k (left), HT (center), and DM (right) estimators, fixing $1/s = 3$ and $\delta = 0.2$ (Uniform network, $x$-axis truncated to $\ell \leq 400$, $y$-axis to $[0, 0.5]$).
  Each curve corresponds to a different fractional parameter: $f_r \in \{0.25, 0.5, 1, 1.5, 2\}$ for Haj\'{e}k/HT and $f_b \in \{0, 0.1, 0.3, 0.5\}$ for DM.
  Haj\'{e}k uniformly dominates HT; DM with $f_b = 0$ is competitive at large $\ell$ but biased at small $\ell$.
  }
  \label{fig:christina_fig3}
\end{figure}

\noindent{\bf Estimator comparison: MSE vs.\ $\ell$.}
\cref{fig:christina_fig3} plots the MSE of the Haj\'{e}k, HT, and DM estimators as a function of $\ell$, fixing $1/s = 3$ and $\delta = 0.2$, with each curve corresponding to a different fractional parameter ($f_r$ for Haj\'{e}k/HT, $f_b$ for DM).

Three observations stand out. First, the Haj\'{e}k estimator uniformly dominates HT across all $\ell$ and $f_r$ values: at $\ell = 100$ and $f_r = 0.5$, the MSE of Haj\'{e}k is roughly an order of magnitude smaller than that of HT, confirming that the self-normalizing denominator provides a large variance-reduction benefit.
Second, the DM estimator with zero burn-in ($f_b = 0$) achieves the lowest MSE among DM variants at large $\ell$ (e.g.\ $\ell \ge 200$), where carry-over effects decay within a single block and the additional data retained by not discarding a burn-in period outweighs the residual bias.
Third, for small $\ell$ (e.g.\ $\ell \le 60$), DM with $f_b = 0$ is substantially biased and its MSE exceeds that of the best HT/Haj\'{e}k variants, while DM with larger $f_b$ reduces bias at the cost of higher
variance, a pattern that mirrors the recency trade-off for HT/Haj\'{e}k.

Taken together, these results reinforce the recommendation to use the
Haj\'{e}k estimator with $f_r \approx 0.5$ as the default choice for spatiotemporal designs on this network.

\noindent{\bf Observations and Discussion.}
The joint spatiotemporal design achieves the lowest MSE across all parameter
configurations, reducing MSE by $58\%$ relative to pure temporal (switchback)
and $97\%$ relative to pure spatial (clustered RD) at their respective optima, confirming the theoretical rate improvement.
The theory establishes that HT is asymptotically superior to DM as
$T\to \infty$; however, in finite samples with $T = 10{,}000$, the best-tuned DM estimator achieves lower MSE than HT for $\ell \ge 80$ (fewer
than 125 blocks).
This is likely because at large $\ell$, there are very few time blocks, and hence the treatment assignment can be more imbalanced (measured the ratio between the size of treated and control units).

\begin{table}[htb]
\centering
\small
\begin{tabular}{lll}
Theoretical prediction & Empirical finding & Status \\
\midrule
Joint $\succ$ pure temporal & Joint MSE $2.4\times$ lower & \checkmark\ Confirmed (\cref{fig:design_viewB}) \\
Joint $\succ$ pure spatial  & Joint MSE $28\times$ lower  & \checkmark\ Confirmed (\cref{fig:design_viewB}) \\
HT $\succ$ DM (asymptotically) & HT $\succ$ DM only at small $\ell$ & $\triangle$\ Finite-sample caveat (\cref{fig:christina_fig3}) \\
MSE U-shaped in $\ell$ & Sweet spot at $\ell \approx 200$--$400$ & \checkmark\ Confirmed (\cref{fig:heatmap_HT,fig:christina_fig3}) \\
Bias--variance trade-off in $r$ & U-shaped MSE vs.\ $r$ & \checkmark\ Confirmed (\cref{fig:christina_fig1}) \\
Haj\'{e}k $\approx$ HT variant & Haj\'{e}k: median $60\times$ lower variance & \checkmark\ Consistent (\cref{fig:christina_fig3}) \\
\end{tabular}
\caption{Theoretical predictions vs.\ empirical findings (Uniform network,
$\delta=0.2$, $T=10{,}000$).}
\label{tab:sim_summary}
\end{table}

This is exactly what Haj\'{e}k adjustment tackles: The ``self-normalizing'' denominator acts as a variance-reduction correction for assignment imbalance.
At $r = 40$, the variance ratio $\mathrm{Var}[\hat\tau_{\mathrm{HT}}] / \mathrm{Var}[\hat\tau_{\mathrm{Haj\acute{e}k}}]$
ranges from $5\times$ to $1028\times$ (median $60\times$), while the bias
ratio remains close to $1$, consistent with the well-known finite-sample superiority of Haj\'{e}k over HT in survey sampling~\citep{hajek1971comment}.
 \cref{tab:sim_summary} summarizes the theoretical predictions and empirical findings.
The practical recommendation is therefore to use Haj\'{e}k, which dominates both HT and DM throughout.

\subsection{Simulation Results: D.C. Hotel Network}
\label{sec:sim_results}

We evaluate the performance of  Horvitz-Thompson (HT)
estimator, Haj\'{e}k, and the Difference-in-Means (DM)
estimators on the network derived from the locations and
capacities of hotels in D.C.
All experiments use a time horizon of $T = 10{,}000$ steps, and the
true GATE is computed via Monte Carlo.

\begin{figure}[htb]
  \centering
  \includegraphics[width=\textwidth]{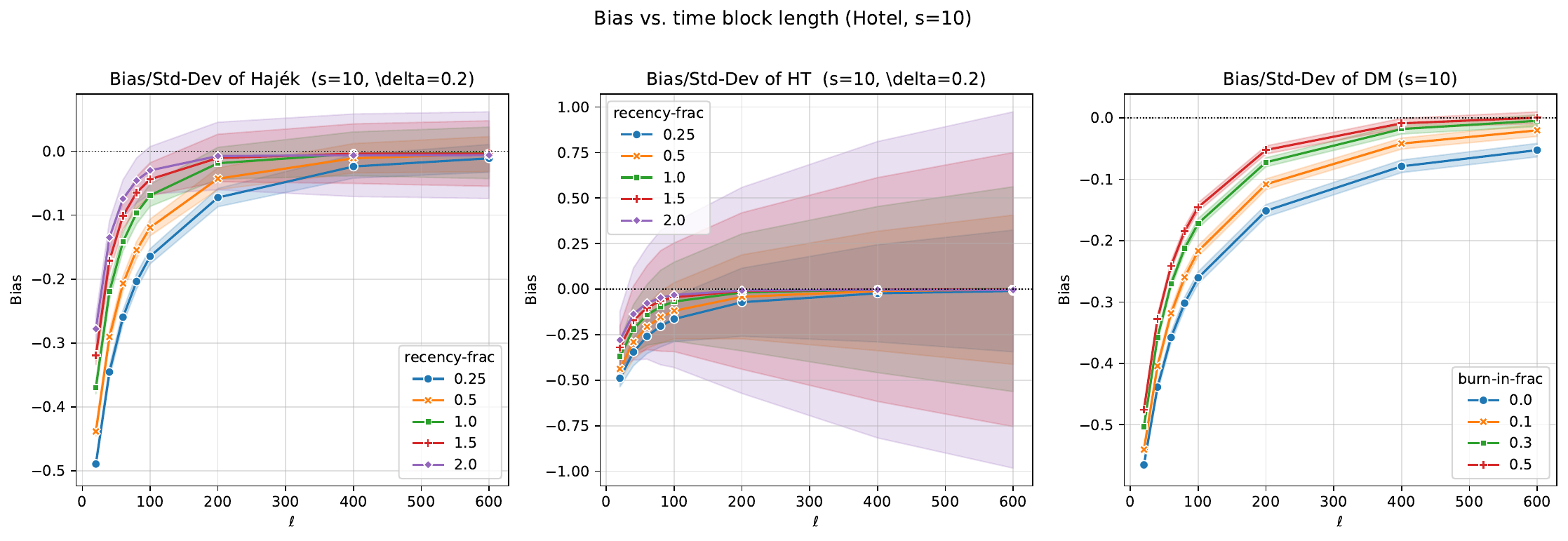}
    \caption{Bias and standard deviation of the Haj\'{e}k, HT, and DM estimators as a function of time-block length $\ell$, fixing $1/s = 10$ and $\delta = 0.2$.
    Each curve corresponds to a different fractional parameter: $f_r \in \{0.25, 0.5, 1, 1.5, 2\}$ for Haj\'{e}k/HT and $f_b \in \{0, 0.1, 0.3, 0.5\}$ for DM. Shaded bands indicate $\pm 1$ standard deviation.}
  \label{fig:ht_heatmap}
\end{figure}

\begin{figure}[htb]
  \centering
  \includegraphics[width=\textwidth]{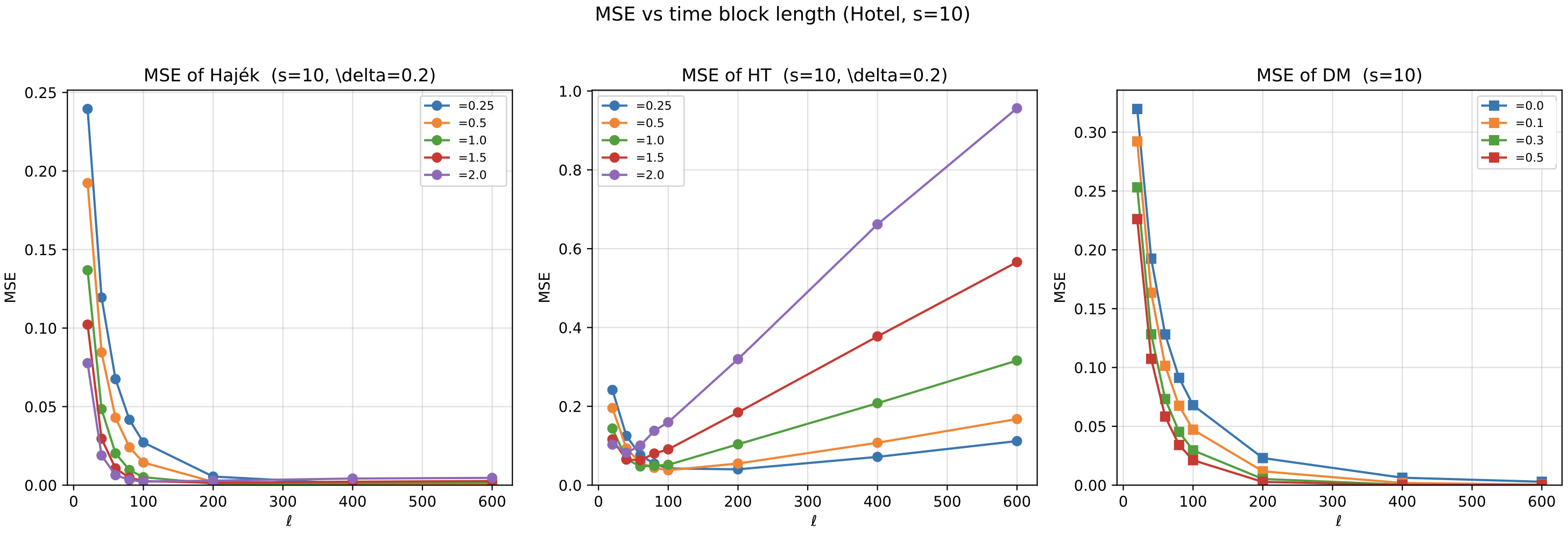}
    \caption{MSE of the Haj\'{e}k, HT, and DM estimators as a function of time-block length $\ell$, fixing $1/s = 10$ and $\delta = 0.2$. Each curve corresponds to a different fractional parameter: $f_r \in \{0.25, 0.5, 1, 1.5, 2\}$ for Haj\'{e}k/HT and $f_b \in \{0, 0.1, 0.3, 0.5\}$ for DM.
    }
    \label{fig:hajek_heatmap}
\end{figure}

\medskip
\noindent{\bf Spatiotemporal Clustered Designs.}
We partition the unit square into a uniform grid of $(1/s) \times (1/s)$ cells, and assign each cell a time block of length $\ell$.
The resulting spatiotemporal clusters are the experimental units of the randomized design.
We vary $1/s \in \{10, 20, 30\}$ and
$\ell \in \{20, 40, 60, 80, 100\}$.
Recall that the HT and Haj\'{e}k estimators are parameterized by a recency $r$ and an exposure-mapping threshold $\delta$; the DM estimator uses a burn-in period $b$.

\cref{fig:ht_heatmap} and \cref{fig:hajek_heatmap} show the bias, standard deviation, and
MSE of the HT estimator as a function of $\ell$ and the recency fraction $r/\ell$, fixing $s= 10$ and $\delta = 0.2$.
We observe that for all three estimators, as the time block length increases, both bias and MSE decrease.
For the Haj\'{e}k and difference-in-means estimators, MSE is driven primarily by bias, whereas for the HT estimator, variance is the dominant component.
We also observe that, the MSEs of HT and DM improve as $\ell$ increases.
This further confirms that the bias is the dominating term.

\begin{figure}[htb]
  \centering
  \includegraphics[width=0.8\textwidth]{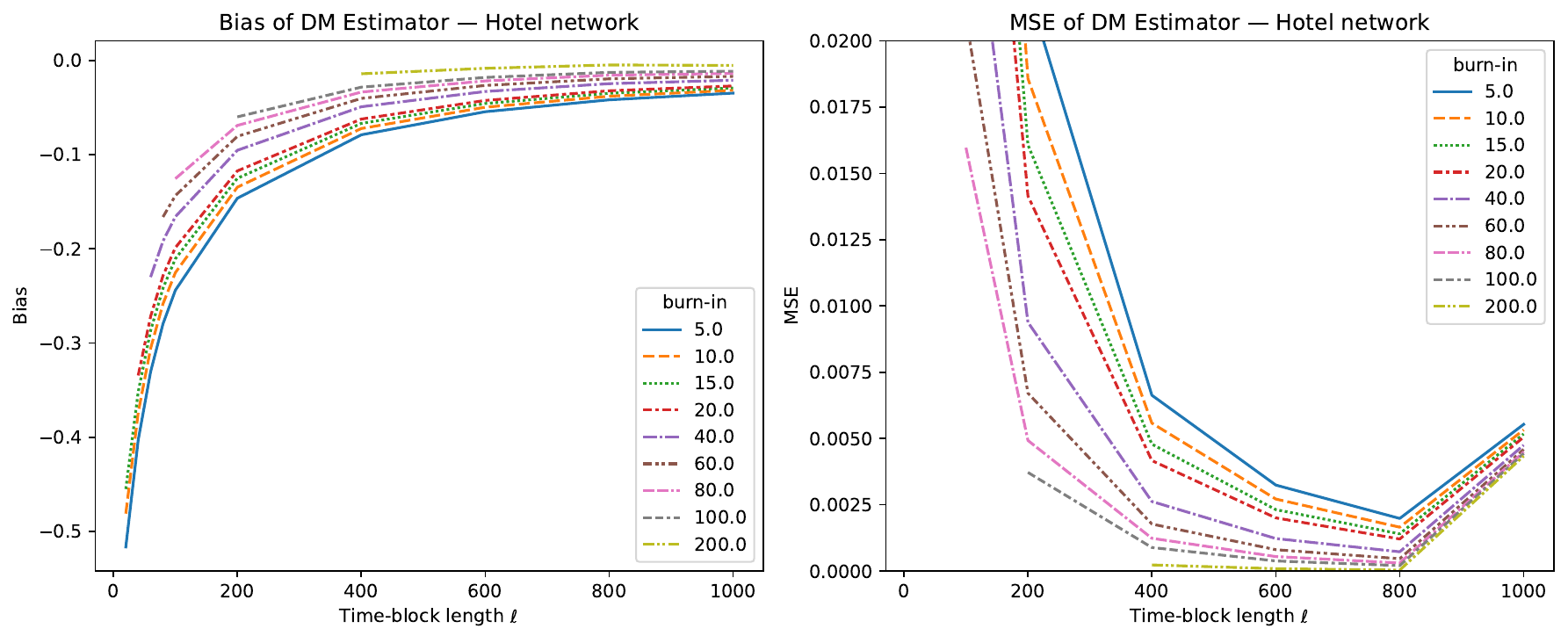}
  \caption{Bias and MSE  of the DM estimator on the Hotel network as a function of $\ell$, for burn-in fractions
    $b/\ell \in \{0, 0.1, 0.3, 0.5\}$.
    Each MSE curve has an interior minimum, and a small positive
    burn-in fraction substantially reduces bias with modest variance
    cost.}
  \label{fig:dm_bias_mse}
\end{figure}

\noindent{\bf Switchback Designs.}
In the switchback setting, the entire network is a single spatial cluster, and vary only
the time-block length $\ell \in \{20, 40, \ldots, 4000\}$.
\cref{fig:dm_bias_mse} shows the bias and MSE of the DM
estimator as a function of $\ell$ for several burn-in fractions
$b \in \{0, 0.1, 0.3, 0.5\} \times \ell$.
Without burn-in ($b = 0$) the estimator is substantially biased for small $\ell$ because carry-over interference from the previous block contaminates the outcome.
Increasing $b$ reduces bias at the cost of discarding observations,
which raises variance and hence MSE for large $\ell$.
The MSE curves therefore exhibit a \emph{sweet spot} in $\ell$: too short a block leaves insufficient washout time; too long a block
wastes data.
A moderate burn-in fraction ($b \approx 0.1\ell$) provides the best MSE.

\begin{figure}[htb]
  \centering
  \includegraphics[width=0.8\textwidth]{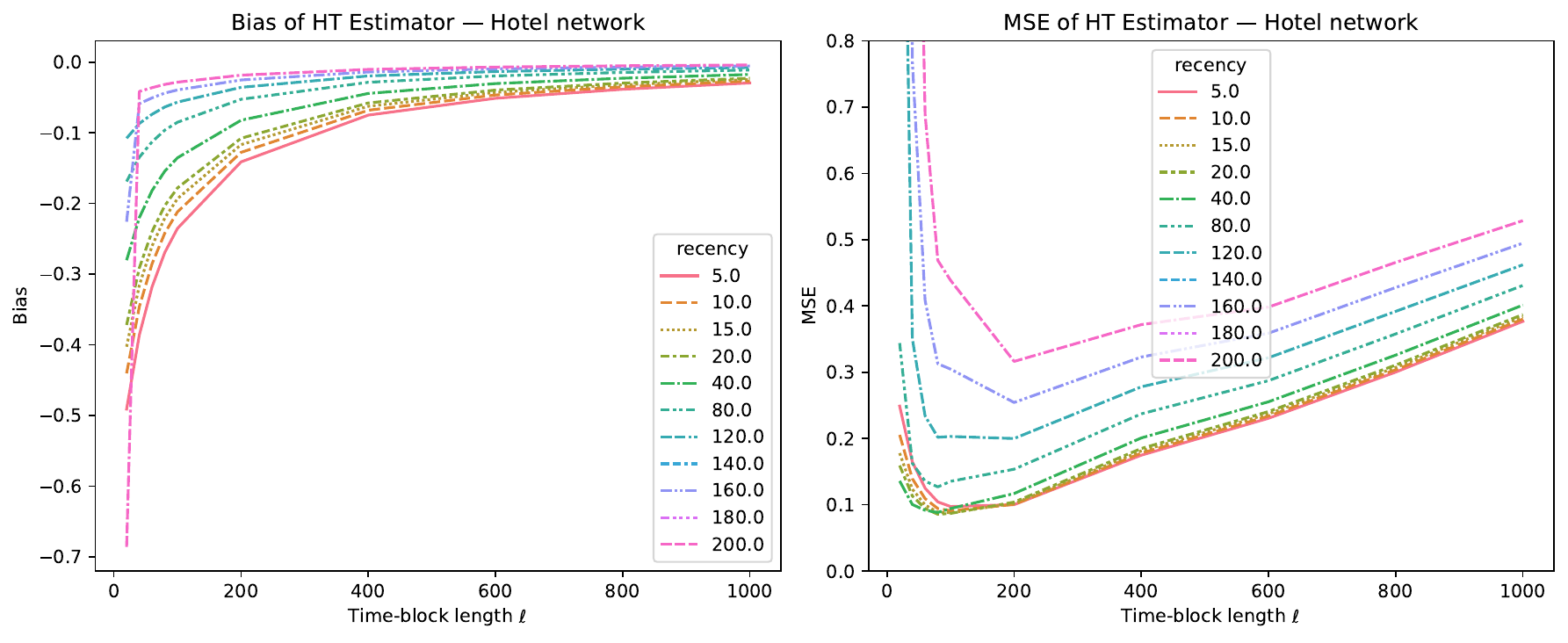}
  \caption{Bias and MSE of the HT estimator as a function of $\ell$, for recency values $r \in \{5, 10, 20, 40, 80, \ldots\}$ with $\delta = 0.2$.
  The HT estimator exhibits large MSE relative to Haj\'{e}k (note the $y$-axis scale), and the optimal $r$ depends on $\ell$.}
  \label{fig:ht_bias_mse_ell}
\end{figure}

\begin{figure}[htb]
  \centering
  \includegraphics[width=0.8\textwidth]{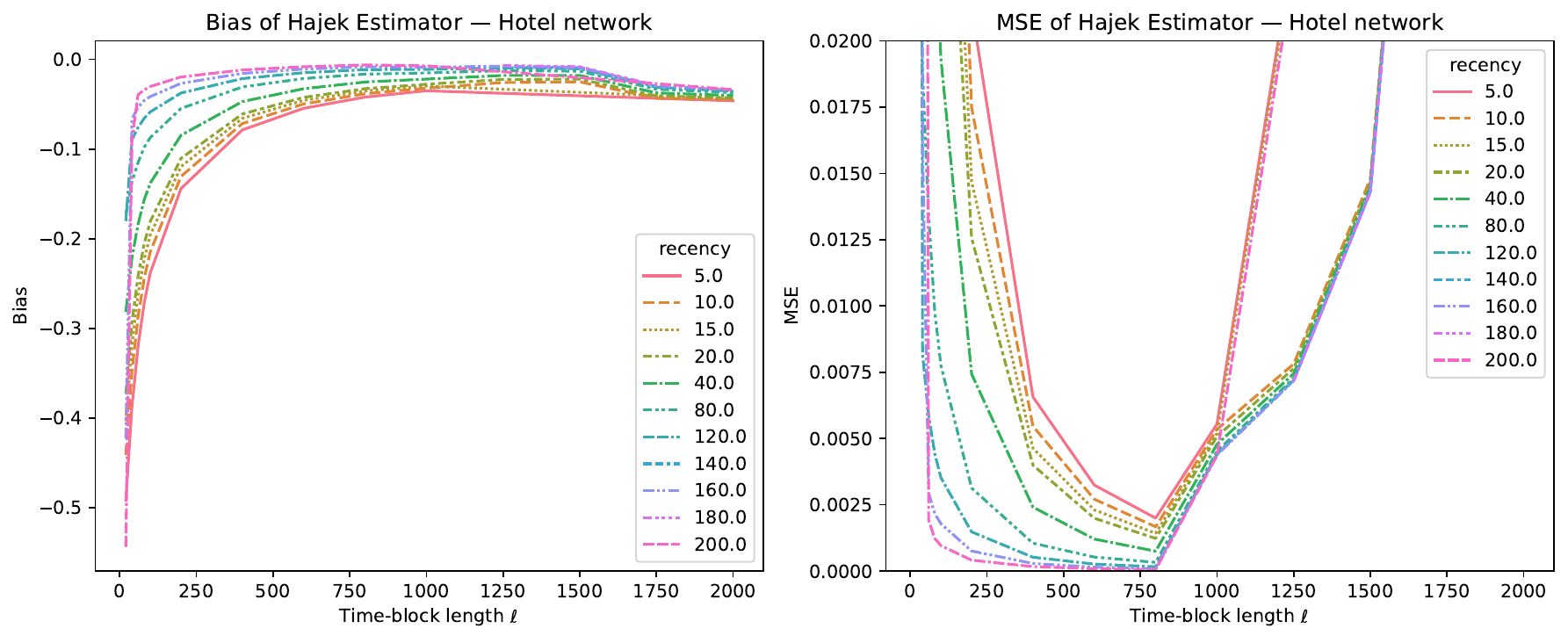}
  \caption{Bias and MSE of the Haj\'{e}k estimator as a function of $\ell$, for recency values $r$ with $\delta = 0.2$.
  The MSE scale ($\le 0.02$) is far smaller than for HT (\cref{fig:ht_bias_mse_ell}), highlighting the efficiency advantage of the Haj\'{e}k correction.}
  \label{fig:hajek_bias_mse_ell}
\end{figure}

\cref{fig:ht_bias_mse_ell} shows the bias and MSE of the HT estimator as a function of $\ell$ for several recency values $r$, fixing $\delta = 0.2$.
The HT estimator is highly sensitive to $r$: small $r$ leads to large bias, while large $r$ inflates variance by assigning near-zero inclusion probabilities to many observations.
The MSE is minimized at an intermediate $r$ that depends on $\ell$, again revealing a sweet spot.
Notably, the MSE of HT is substantially larger than that of Haj\'{e}k (compare the $y$-axis scales: HT reaches $0.8$ whereas Haj\'{e}k
stays below $0.02$), confirming that the self-normalization in Haj\'{e}k is critical.

\cref{fig:hajek_bias_mse_ell} shows the bias and MSE of the
Haj\'{e}k estimator as a function of $\ell$.
The MSE is uniformly below $0.02$—roughly $40$ times lower than
the HT estimator at comparable parameter settings, demonstrating the
substantial efficiency gain from self-normalization.
The bias behavior is qualitatively similar to HT, but the variance
reduction is dramatic.

\medskip
\noindent{\bf Summary of Insights.}
The simulations study leads to three main takeaways.
\begin{enumerate}
  \item  \textbf{MSE: Haj\'{e}k $\ll$ HT $<$ DM.}
    The Haj\'{e}k estimator consistently achieves the lowest MSE,
    often by an order of magnitude relative to HT and DM.
    The self-normalizing denominator effectively controls the
    variance inflation caused by near-zero inclusion probabilities.
  \item \textbf{MSE Sweet Spots.}
    Every estimator exhibits an interior optimum in its design
    parameters $(\ell, r, b)$.  The optimal $\ell$ balances
    interference washout (favoring large $\ell$) against sample
    efficiency (favoring small $\ell$).  The optimal recency $r$
    balances bias (favoring large $r$) against variance (favoring
    small $r$), and scales approximately linearly with $\ell$.
  \item \textbf{Diminishing Gains of Having More Spatial clusters.}
    Comparing the spatiotemporal and pure switchback settings, adding
    spatial clusters ($1/s > 1$) reduces MSE further, but the
    marginal benefit diminishes as $1/s$ grows beyond $\sim 20$.
\end{enumerate}

\section*{Acknowledgments}
The authors thank the reviewers for their valuable comments, which substantially improved the manuscript. During the revision, the authors used Claude to  check the manuscript for typographical and notational consistency and to assist with editing. All suggestions were reviewed by the authors, who take full responsibility for the content. No AI tool was used to generate the theoretical or simulation results.

%% file: sections/apdx.tex
\section{Bias Analysis}\label{sec:bias_proof}

Next, we show that \cref{assu:mixing} implies a bound on how the outcome distribution under $\mathbb P[\cdot \mid W=w]$ can vary in $\bf w$.

\begin{lemma}[Decaying Temporal Interference]
\label{prop:TVdecay}
Consider any $t,m\in [T]$ with $1\le m<t\le T$ and $i\in U$.
Suppose $w, w'\in \{0,1\}^{N\times T}$ are identical on $\N(i) \times [t-m, t]$. Then,
\[
    d_{\rm TV} \lb(\mathbb P[Y_{it}\in\cdot\mid W=w], \ \mathbb P[Y_{it}\in\cdot\mid W=w']\rb)
    \le e^{-m/t_{\rm mix}}.
\]
\end{lemma}
\proof For any $i,t$, we denote $f_{it} = \mathbb P[S_{it} \in \cdot \mid W=w]$ and $f'_{it}= \mathbb P[S_{it}\in\cdot\mid W=w']$.
Then, for any $s\in [T]$, by the Chapman–Kolmogorov equation,
\[
    f_{i,s+1} = f_{is} P^{w_{\N(i)},s}_{is}\quad \text{and} \quad f'_{i,s+1} = f'_{is} P^{w'_{\N(i)},s}_{is}.
\]
As $w,w'$ are identical on $\N(i) \times [t-m, t]$, it follows that for $t-m\le s\le t$,
\begin{align*}
    d_{\rm TV}(f_{i, s+1},\ f'_{i, s+1})
    &= d_{\rm TV}\lb(f_{is} P^{w_{\N(i)},s}_{is},\ f'_{is} P^{w'_{\N(i)},s}_{is}\rb)\\
    &= d_{\rm TV}\lb(f_{is} P^{w_{\N(i)},s}_{is},\ f'_{is} P^{w_{\N(i)},s}_{is}\rb)\\
    &\le e^{-1/t_{\rm mix}} \cdot d_{\rm TV}\lb(f_{is}, f'_{is}\rb),
\end{align*}
where we used $w'_{\N(i),s}=w_{\N(i),s}$ in the second equality and \cref{assu:mixing} in the inequality.
Applying the above for all $s=t-m,\dots,t$, it follows that
\begin{align*}
    d_{\rm TV} \lb(f_{it},f'_{it}\rb)
    & \le e^{-m/ t_{\rm mix}} \cdot d_{\rm TV}(f_{i,t-m},f'_{i,t-m})
    \le e^{-m/ t_{\rm mix}}.
\end{align*}
Finally, by the data processing inequality, it follows that \begin{align*}
    d_{\rm TV} \lb(\mathbb P[Y_{it}\in\cdot\mid W=w],\ \mathbb P[Y_{it}\in\cdot\mid W=w'] \rb)
    &\le d_{\rm TV} \lb(f_{it},f'_{it}\rb) \le e^{-m/ t_{\rm mix}}.
\end{align*}
\qed

\cref{prop:TVdecay} leads to the following bound on the bias.
\begin{lemma}[Per-unit Bias]\label{lem:bias_per_round} For any $a\in \{0,1\}$, $r>0$, $i\in U$ and $t\in [T]$, we have \[\lb|\ho{E}\lb[\wh Y_{ita}^r\rb] - \ho{E}\lb[Y_{it}\mid  W = a \cdot {\bf 1}\rb] \rb| \le e^{-r /t_{\rm mix}}.\]
\end{lemma}
\proof  For any $a\in \{0,1\}$, $i\in U$ and $t\in [T]$, we have
\begin{align*}
\ho{E}\lb[\wh Y_{ita}^r\rb] &= \ho{E}\lb[\frac{X_{ita}^r}{p_{ita}^r} Y_{it} \middle| X_{ita}^r = 1\rb] \ho{P}\lb[X_{ita}^r = 1\rb] + \ho{E}\lb[\frac{X_{ita}^r}{p_{ita}^r} Y_{it} \middle| X_{ita}^r = 0\rb] \ho{P}\lb[X_{ita}^r = 0\rb] \\
&= \ho{E}\lb[\frac{X_{ita}^r}{p_{ita}^r} Y_{it} \middle| X_{ita}^r = 1\rb] p_{ita}^r +0 \\
& = \ho{E} \lb[Y_{it} \mid X_{ita}^r =1\rb].
\end{align*}
Note that $X_{ita}^r=1$ implies that $w=a\cdot{\bf 1}$ on  $\N(i)\times [t-r,t]$.
Therefore, by \cref{prop:TVdecay} with $m=r$, we obtain
\[
    d_{\rm TV} \lb(\mathbb P[Y_{it}\in\cdot\mid W=w],\ \mathbb P[Y_{it}\in\cdot\mid W={a\cdot \bf 1}]\rb)
    \le e^{-r/t_{\rm mix}}.
\]
Because $\mu_{it} \in [0,1]$ it follows that
\[\lb|\ho{E}\lb[\wh Y_{ita}^r\rb] - \ho{E}\lb[Y_{it}\mid W =a\cdot {\bf 1}\rb]\rb| = \lb|\ho{E} \lb[Y_{it} \mid X_{ita}^r=1\rb] - \ho{E}\lb[Y_{it}\mid W = a\cdot {\bf 1}\rb] \rb| \le e^{-r/t_{\rm mix}}.\eqno \qed\]

Now we are prepared to upper bound the bias.
Recall that $\Delta = \frac 1 {NT} \sum_{i,t} \Delta_{it}$ and $\wh\Delta^r = \frac 1 {NT} \sum_{i,t} \wh\Delta_{it}^r$.

\medskip
{\bf \noindent Proof of \cref{prop:bias}.}
By \cref{lem:bias_per_round},
\begin{align*}
\lb|\ho{E} \lb[\wh\Delta^r\rb] - \Delta\rb| & \le \frac 1{NT}\sum_{i,t} \lb|\Delta_{it} - \ho{E}\lb[\wh \Delta_{it}^r\rb]\rb| \\
& \le \frac 1{NT} \sum_{a\in \{0,1\}}\sum_{i,t}\lb|\ho{E}\lb[\wh Y_{ita}^r\rb] - \ho{E}\lb[Y_{it}\mid W=a\cdot {\bf 1}\rb] \rb|\\
&\le 2e^{-r/t_{\rm mix}}.\tag*\qed
\end{align*}

\section{Variance Analysis} \label{sec:variance_proof}

We start with a bound that holds for all pairs of units.
Note that $p_{it0}^r = p_{it1}^r$ for any $i,t,r$, since treatment and control are assigned with equal probabilities. We will thus suppress $a$.
\begin{lemma}[Covariance Bound]
\label{lem:close-by}
For any $r,\ell\ge 0$, $i,i'\in U$ and $t,t'\in [T]$, we have
\[
    {\rm Cov}\lb(\wh\Delta_{it}^r, \wh\Delta_{i't'}^r\rb) \le \frac {4(1+\sigma^2)}{ \sqrt{p_{it}^r\cdot p_{i't'}^r}}.
\]
\end{lemma}
\proof
Expanding the definition of $\wh \Delta_{it}^r$, we have
\begin{align}
\label{eqn:121723}
    {\rm Cov}\lb(\wh\Delta_{it}^r, \wh\Delta_{i't'}^r\rb)
    & = {\rm Cov}\lb(\frac{X_{it1}^r}{p_{it1}^r}Y_{it} - \frac{X_{it0}^r}{p_{it0}^r}Y_{it},\ \frac{X_{it1}^r}{p_{i't'1}^r}Y_{i't'}- \frac{X_{i't'0}^r}{p_{i't'0}^r}Y_{i't'}\rb)\notag \\
    &\le \sum_{a,a'\in \{0,1\}} \lb|{\rm Cov}\lb(\frac{X_{ita}^r}{p_{ita}^r}Y_{it},\ \frac{X_{i't'a'}^r}{p_{i't'a'}^r}Y_{i't'}\rb)\rb| \notag\\
    &\le \sum_{a,a'\in \{0,1\}} \frac 1{p_{ita}^r}\frac 1{p_{i't'a'}^r} \lb|{\rm Cov}\lb(X_{ita}^r Y_{it},\ X_{i't'a'}^r Y_{i't'}\rb)\rb|  \notag\\
    &\le \sum_{a,a'\in \{0,1\}} \frac 1{p_{ita}^r} \frac 1{p_{i't'a'}^r} \sqrt {{\rm Var}[X_{ita}^r Y_{it}]} \sqrt {{\rm Var}[X_{i't'a'}^r Y_{i't'}]},
\end{align}
where the last inequality is by the Cauchy-Schwarz inequality.
Note that $X_{ita}^r$ is binary and the outcome function $\mu$ takes values in [0,1], such that

\begin{align*}
\rm Var[X_{ita}^r Y_{it}] &=
{\rm Var}[ \ho{E} [X_{ita}^r Y_{it} ~|~ W, S_{it}] ]
+ \ho{E} [{\rm Var}[X_{ita}^r Y_{it}~|~ W, S_{it}]  ] \\
&\le {\rm Var}[X_{ita}^r \mu_{it}\lb(S_{it}, W_{\N(i),t})\rb)] + \ho{E} [X_{ita}^r \sigma^2] \\
&\le \ho{E}[X_{ita}^r \mu_{it}\lb(S_{it},(W_{\N(i),t})\rb)^2]
+ p_{ita}^r \sigma^2\\
&\le p_{ita}^r (1 + \sigma^2).
\end{align*}
It follows that
\begin{align*}
    {\rm Cov}\lb(\wh\Delta_{it}^r, \wh\Delta_{i't'}^r\rb)
    \le \lb(\frac 1{\sqrt{p_{it0}^r}} + \frac 1{\sqrt{p_{it1}^r}} \rb) \lb(\frac 1{\sqrt{p_{i't'0}^r}} + \frac 1{\sqrt{p_{i't'1}^r}}\rb)(1+\sigma^2)
    = \frac{4(1+\sigma^2)}{\sqrt{p_{it}^r\cdot p_{i't'}^r}}.\tag*\qed
\end{align*}

The above bound alone is not sufficient for our analysis, as it does not take advantage of the rapid mixing property.
Next, in \cref{lem:cov_outcomes}, we show that for any units that are far apart in time, the covariance of their HT terms decays {\bf exponentially} in their temporal distance. In order to show that, we first show that if the realization of one random variable has little impact on the (conditional) distribution of another random variable, then they have low covariance.

\begin{lemma}[Low Interference in Conditional Distribution Implies Low Covariance]\label{lem:low_cov}
Let $U,V$ be two random variables and $g,h$ be real-valued functions defined on their respective realization spaces.
If for some $\delta>0$, we have \[d_{\rm TV}(\mathbb P[U\in\cdot\mid V],\ \mathbb P[U\in\cdot])\le \delta \quad V\text{-almost surely },\]
then,
\[
    |{\rm Cov}(g(U),h(V))| \le 2 \delta \cdot \|h(V)\|_1 \cdot  \|g(U)\|_\infty
\]
\end{lemma}
\proof
Denote by $\mu_{U,V},\mu_U,\mu_V,\mu_{U\mid V=v}$ the probability measures of $(U,V)$, $U$, $u$, and $U$ conditioned on $V=v$, respectively. We then have
\begin{align*}
|{\rm Cov}(g(U),h(V))| &= |\ho{E}\lb [g(U)  h(V)\rb] - \ho{E}[g(U)] \ho{E} [h(V)]| \notag\\
& =\lb|\int_v h(v) \lb(\int_u g(u) \lb(\mu_{U\mid V=v}(du) - \mu_U(du) \rb)\rb)\mu_V(dv)\rb|\\
&\le \int_v |h(v)| \cdot \|g(U)\|_\infty \cdot 2 d_{\rm TV}(\mathbb P(U\in\cdot\mid V=v), \mathbb P(U\in\cdot))\mu_V(dv)\\
&\le 2\|h(V)\|_1\cdot \|g(U)\|_\infty \cdot \delta.\tag*\qed
\end{align*}

Viewing $V,U$ as outcomes in different rounds, we use the above to bound the covariance in the outcomes in terms of their temporal distance.

\begin{lemma}[Covariance of Outcomes] \label{lem:cov_outcomes}
For any $A \sse \{0,1\}^{N\times T}$, $i,i'\in U$ and $t,t'\in [T]$, we have
\[
    |{\rm Cov} (Y_{it}, Y_{i't'}\mid W\in A)| \le e^{-{|t-t'|/{t_{{\rm mix}}}}}.
\]
\end{lemma}
\proof
Wlog assume $t'<t$.
Recall that $\epsilon_{it} =Y_{it}-\mu_{it}(S_{it},W_{\N(i),t})$, so
\begin{align*}
    {\rm Cov}(Y_{i't'},Y_{it}\mid W\in A)
    & = {\rm Cov}(\mu_{i't'}(S_{i't'}, W_{i't'}),\mu_{it}(S_{it}, W_{it}) \mid W\in A)\\
    &\quad + {\rm Cov}(\mu_{i't'}(S_{i't'}, W_{i't'}),\epsilon_t\mid W \in A) \\
    & \quad +{\rm Cov}(\epsilon_{i't'},\mu_{it}(S_{it},W_{it})\mid W\in A) \\
    &\quad + {\rm Cov}(\epsilon_{i't'},\epsilon_{it}\mid W\in A).
\end{align*}
The latter three terms are zero by the exogenous noise assumption (in terms of covariances).
By \cref{prop:TVdecay} and the triangle inequality for $d_{\rm TV}$, for any $s\in\mathcal{S}$,  $w\in A$, we have
\begin{align}
\label{eqn:121423b}
    d_{\rm TV}\lb(\ho{P} \lb[(S_{it},W_{\N(i), t})\in\cdot \mid W\in A \rb],\ \ho{P} \lb[(S_{it},W_{\N(i),t})\in \cdot \mid (S_{i't'},W_{\N(i'),t'})=(s,w)\rb]\rb)
    \le e^{-{(t-t')/{t_{{\rm mix}}}}},
\end{align}
We obtain the desired inequality from \cref{eqn:121423b} by applying    \cref{lem:low_cov}, with
$(S_{it},W_{it})$ in the role of ``$U$'', $(S_{i't'},W_{i't'})$ in the role of ``$V$'', $\mu_{it},\mu_{i't'}$ in the roles of ``$g$'' and ``$h$'', $e^{-(t-t')/t_{\rm mix}}$ in the role of ``$\delta$'',
and observing that $\|\mu_{it}\|_\infty, \|\mu_{i't'}\|_\infty\le 1$ (and hence $\|\mu_{i't'}(V)\|_1\le 1$ for any random variable $V$). \qed

Using the result in \cref{lem:low_cov} that the outcomes have low covariance if they are far apart in time, we show that the covariance between the HT terms, $\wh \Delta^r_{it}$ is also low. Each of the HT terms are a product of the outcome and the exposure mapping.

\begin{lemma}[Covariance of Far-apart HT terms]
\label{lem:cov_bound}
Suppose $i,i'\in U$ and $t,t'\in [T]$ satisfy $\lceil t'/\ell\rceil > \lceil t/\ell\rceil$ and $t+\ell < t'-r$,
then
\[
    {\rm Cov}\lb(\wh\Delta_{it}^r, \wh\Delta_{i't'}^r\rb) \le 8e^{-|t'-t|/t_{\rm mix}}.
\]
\end{lemma}
\proof
By the definition of  covariance,
\begin{align}\label{eqn:121423}
&\quad \ {\rm Cov}\lb(X_{ita}^r Y_{it},\ X_{i't'a'}^r Y_{i't'}\rb) \notag\\
&= \ho{E}[X_{ita}^r Y_{it}\cdot X_{i't'a'}^r Y_{i't'}]  - \ho{E}[X_{ita}^r Y_{it}]\cdot \ho{E}[X_{i't'a'}^r Y_{i't'}]\notag\\
&= \ho{E}[Y_{it} Y_{i't'}\mid X_{ita}^r = X_{i't'a'}^r =1]\cdot \ho{P}[X_{ita}^r = X_{i't'a'}^r =1] \notag \\
& \quad \quad - \ho{E}[Y_{it} \mid X_{ita}^r=1]
\cdot \ho{P}[X_{ita}^r=1] \cdot \ho{E}[Y_{i't'}\mid X_{i't'a'}^r=1] \cdot
\ho{P}[X_{i't'a'}^r =1].
\end{align}
To further simplify, observe that for any (possibly identical) $a,a'\in \{0,1\}$, since $t,t'$ lie in distinct time blocks and are more than $2r$ apart, we have
\[
    X_{ita}^r\ind X_{i't'a'}^r.
\]
So,
\begin{align}\label{eqn:052525z}
\eqref{eqn:121423} &= \ho{E}[Y_{it} Y_{i't'}\mid X_{ita}^r = X_{i't'a'}^r =1]\cdot \ho{P}[X_{ita}^r=1] \cdot
\ho{P}[X_{i't'a'}^r =1] \notag \\
&\quad\quad - \ho{E}[Y_{it} \mid X_{ita}^r=1]
\cdot \ho{E}[Y_{i't'}\mid X_{i't'a'}^r=1]
\cdot \ho{P}[X_{ita}^r=1] \cdot
\ho{P}[X_{i't'a'}^r =1]  \notag\\
& = \ho{P}[X_{ita}^r=1] \cdot
\ho{P}[X_{i't'a'}^r =1] \cdot \lb({\rm Cov}(Y_{it}, Y_{i't'}\mid X_{ita}^r = X_{i't'a'}^r =1) + \epsilon+ \epsilon' + \epsilon\epsilon'\rb),
\end{align}
where
\[
    \epsilon := \ho{E}[Y_{it} \mid X_{ita}^r=1]  - \ho{E}[Y_{it} \mid X_{ita}^r = X_{i't'a'}^r=1]
\]
and
\[
    \epsilon' := \ho{E}[Y_{i't'} \mid X_{i't'a'}^r=1] -\ho{E}[Y_{i't'} \mid X_{ita}^r = X_{i't'a'}^r=1].
\]
Equivalently, \cref{eqn:052525z} can be rewritten as
\begin{align}\label{eqn:052425c}
    {\rm Cov}\lb(\frac{X_{ita}^r Y_{it}}{p_{ita}^r},\ \frac{X_{i't'a'}^r Y_{i't'}}{p_{i't'a'}^r}\rb)
    = {\rm Cov}(Y_{it}, Y_{i't'}\mid X_{ita}^r = X_{i't'a'}^r =1) + \epsilon+ \epsilon'+\epsilon\epsilon'
\end{align}
We split the rest of the analysis into three parts.
\begin{enumerate}
\item Note that the state at any time only depends on the treatment assignment {\bf in the  past}  and that $t<t'$, so
\[Y_{it}\ind X_{i't'a'}^r \mid X_{ita}^r,\]
and hence \begin{align}\label{eqn:052425a} \epsilon = 0.
\end{align}
\item Furthermore, by \cref{prop:TVdecay}, we have
\[
    d_{\rm TV}(\ho{P}[Y_{i't'}\in \cdot \mid X_{ita}^r=1],\ \ho{P}[Y_{i't'}\in \cdot \mid X_{ita}^r=X_{i't'a'}^r=1])\le e^{-|t-t'|/t_{\rm mix}},
\]
and therefore
\begin{align}\label{eqn:052425b}
|\epsilon'| \le e^{-|t-t'|/t_{\rm mix}}.
\end{align}
\item By \cref{lem:cov_outcomes}, we have \begin{align}\label{eqn:052525y} \left|{\rm Cov}\lb(\frac{X_{ita}^r Y_{it}}{p_{ita}^r},\ \frac{X_{i't'a'}^r Y_{i't'}}{p_{i't'a'}^r}\rb)\right| \le e^{-|t-t'|/t_{\rm mix}}.\end{align}
\end{enumerate}
Combining, \cref{eqn:052425a,eqn:052425b,eqn:052525y}, we obtain
\[\eqref{eqn:052425c} \le 2e^{-|t-t'|/t_{\rm mix}}.\]
Therefore, by the triangle inequality, we conclude that
\[{\rm Cov}\lb(\wh \Delta_{it}^r,\wh \Delta_{i't'}^r \rb) \le \sum_{(a,a')\in \{0,1\}^2} \lb|{\rm Cov}\lb(\frac{X_{ita}^r Y_{it}}{p_{ita}^r},\ \frac{X_{i't'a'}^r Y_{i't'}}{p_{i't'a'}^r}\rb) \rb| \le 8e^{-|t-t'|/t_{\rm mix}}.\eqno\qed\]

\begin{remark} The restriction that $t,t'$ are both farther than $2r$ apart in time \emph{and} lies in distinct time-blocks is both necessary for the above exponential covariance bound.
To illustrate, fix an $i\in [N]$.
Suppose intervals $[t-r,t], [t'-r,t]$ lie in the same  time block denoted $\cal T$.
Then, the exposure mappings $X_{ita}^r= X_{it'a}^r:=X$ a.s..
Consequently, the HT terms for $it$ and $it'$ are both affected by $X$, and their covariance may {\bf not} decay exponentially in $r$.

More concretely,  it is not hard to construct an example where $Y_{i,\cal T}\equiv 1$ if the time block $\cal T$ is assigned treatment $1$, and $Y_{i,\cal T}\equiv 0$ otherwise.
In this case, both $\wh \Delta_{it}^r= \wh \Delta_{it'}^r$ a.s. and they follow ${\rm Ber}(\frac 12)$.
Thus,
\[
    {\rm Cov}\lb(\wh \Delta_{it}^r,\wh \Delta_{i't'}^r\rb) = \frac 14.
\]
\end{remark}

Finally we put all the pieces together to obtain the final variance bound.

\subsubsection*{Proof of \cref{prop:var}.}
Write
\begin{align}\label{eqn:062324}
{\rm Var}\lb(\wh\Delta^r\rb) = {\rm Var}\lb(\frac 1{NT} \sum_{(i,t)\in U\times [T]} \wh\Delta_{it}^r \rb) = \frac 1{N^2 T^2} \sum_{i,t} \sum_{i',t'}
{\rm Cov}\lb(\wh \Delta_{i't'}^r, \wh \Delta_{it}^r\rb).
\end{align}
We need two observations to decompose the above.
First, by the definition of the dependence graph, if $i\not \ind i'$, then
\[
{\rm Cov}\lb(\wh \Delta_{i't'}^r,\ \wh \Delta_{it}^r\rb) = 0.
\]
Thus, for each $i$, we may consider only the units $i'$ with $i\not \ind i'$.
Second, observe that if $|t-t'|> r+\ell$, then we can apply \cref{lem:cov_bound} to bound the covariance term exponentially in $|t-t'|$.
Combining, we can decompose \cref{eqn:062324} into close-by (``C'') and far-apart (``F'') pairs as
\begin{align}
\eqref{eqn:062324} = \frac 1{N^2 T^2} \sum_{(i,t)} \lb(C_{it}+F_{it}\rb)
\end{align} where
\[C_{it} := \sum_{t':|t'-t|\le \ell+r} \sum_{i': i\not \ind i'}
{\rm Cov}\lb(\wh\Delta_{it}^r, \wh\Delta_{i't'}^r \rb)\] and
\[F_{it} := \sum_{t':|t'-t|> \ell+r} \sum_{i': i\not \ind i'} {\rm Cov}\lb(\wh\Delta_{it}^r, \wh\Delta_{i't'}^r \rb).\]
To further analyze the above, fix any $(i,t)\in U\times [T]$.

\noindent{\bf Part I: Bounding $C_{it}$.}
By \cref{lem:close-by},
\begin{align*}
C_{it}
&\le \sum_{t':|t'-t|\le \ell+r} \sum_{i': i\not \ind i'}  \frac {4(1+\sigma^2)}{\sqrt{p_{it}^r\cdot p_{i't'}^r}}\\
& \le (\ell + r) \sum_{i': i\not \ind i'} \frac{4(1+\sigma^2)}{ \sqrt{p_i^{\rm min}\cdot p_{i'}^{\rm min}}}.
\end{align*}
Summing over all $(i,t)$, we have
\begin{align}
\label{eqn:121423d}
\sum_{i,t} C_{it}
& \le \sum_{i,t} (r+\ell) \sum_{i':i\not \ind i'} {\frac{4(1+\sigma^2)}{\sqrt{p_i^{\rm min}\cdot p_{i'}^{\rm min}}}} \notag\\
&\le T (r+\ell)\sum_{(i,i'):i\not \ind i'} \frac{4(1+\sigma^2)}{\sqrt{p_i^{\rm min}\cdot p_{i'}^{\rm min}}}.
\end{align}

\noindent{\bf Part II: Bounding $F_{it}$.}
By \cref{lem:cov_bound},
\begin{align}\label{eqn:121423e}
F_{it} &=\sum_{i': i\not \ind i'} \sum_{t':|t'-t|> \ell +r} {\rm Cov} \lb(\wh\Delta_{it}^r,\wh\Delta_{i't'}^r\rb)\notag\\
&\le\sum_{i': i\not \ind i'} 2\int_{\ell +r}^\infty 4 e^{-z/t_{\rm mix}}\ dz \notag\\
&=8t_{\rm mix} e^{-(r+\ell)/t_{\rm mix}}\cdot {|\{i': i\not\ind i'\}|},
\end{align}
where the ``2'' in the inequality arises since $s$ may be either greater or smaller than $t$.
Combining \cref{eqn:121423d,eqn:121423e}, we conclude that
\begin{align*}
    {\rm Var}\lb(\wh \Delta^r\rb) & \le \frac 1{N^2 T^2}\lb(\sum_{i,t} C_{it} + \sum_{i,t} F_{it}\rb)\\
    & \le \frac 1{N^2 T^2} \lb({4 (r + \ell)} T \sum_{(i,i'): i\not \ind i'} \frac {1+\sigma^2}{\sqrt{p_i^{\rm min} p_{i'}^{\rm min}}} + 8T t_{\rm mix} e^{-\frac {(r+\ell)}{t_{\rm mix}}} \left|\{(i,i'): i\not \ind i'\}\right| \rb).
\end{align*}
Note that $p_i^{\rm min}$ must be upper bounded by $1/2$, such that $1/\sqrt{p_i^{\rm min} p_{i'}^{\rm min}} \geq 2$. For $r \geq t_{\rm mix}$, the second term is dominated by the first so that
\begin{align*}
    {\rm Var}\lb(\wh \Delta^r\rb) & \le 8(1+\sigma^2)\lb(\frac{r+ \ell}{T} \rb) \lb(\frac{1}{N^2} \sum_{(i,i'): i\not \ind i'} \frac{1}{ \sqrt{p_i^{\rm min} p_{i'}^{\rm min}}}\rb).\tag*\qed
\end{align*}

%% file: references.bib
@article{missault2025robust,
  title={Robust and efficient multiple-unit switchback experimentation},
  author={Missault, Paul and Masoero, Lorenzo and Delb{\'e}, Christian and Richardson, Thomas and Imbens, Guido},
  journal={arXiv preprint arXiv:2506.12654},
  year={2025}
}

@article{hajek1971comment,
  title={Comment on “an essay on the logical foundations of survey sampling, part one”},
  author={H{\'a}jek, Jaroslav},
  journal={The foundations of survey sampling},
  volume={236},
  year={1971},
  publisher={Holt, Rinehart, and Winston Toronto}
}

@misc{dc_hotels_2025,
  author = {{District of Columbia Government}},
  title = {Hotels},
  year = {2025},
  publisher = {ArcGIS Hub},
  url = {https://hub.arcgis.com/datasets/dcgis::hotels/about},
  note = {Accessed: 2026-02-17}
}

@article{kandiros2024conflict,
title={The Conflict Graph Design: Estimating Causal Effects under Arbitrary Neighborhood Interference},
  author={Kandiros, Vardis and Pipis, Charilaos and Daskalakis, Constantinos and Harshaw, Christopher},
  journal={arXiv preprint arXiv:2411.10908},
  year={2024}
}

@article{ugander2023randomized,
  title={Randomized graph cluster randomization},
  author={Ugander, Johan and Yin, Hao},
  journal={Journal of Causal Inference},
  volume={11},
  number={1},
  pages={20220014},
  year={2023},
  publisher={De Gruyter}
}

@article{wang2021causal,
  title={Causal Inference with Panel Data under Temporal and Spatial Interference},
  author={Wang, Ye},
  journal={arXiv preprint arXiv:2106.15074},
  year={2021}
}

@article{sneider19,
  title={Experiment Rigor for Switchback Experiment Analysis},
  author={Sneider, Carla Sneider and Tang, Yixin},
  url={https://careersatdoordash.com/blog/experiment-rigor-for-switchback-experiment-analysis/
},
  year={2019}
}

@article{savje2024causal,
  title={Causal inference with misspecified exposure mappings: separating definitions and assumptions},
  author={S{\"a}vje, Fredrik},
  journal={Biometrika},
  volume={111},
  number={1},
  pages={1--15},
  year={2024},
  publisher={Oxford University Press}
}

@article{li2023experimenting,
  title={Experimenting under stochastic congestion},
  author={Li, Shuangning and Johari, Ramesh and Kuang, Xu and Wager, Stefan},
  journal={arXiv preprint arXiv:2302.12093},
  year={2023}
}

@article{li2022network,
  title={Network interference in micro-randomized trials},
  author={Li, Shuangning and Wager, Stefan},
  journal={arXiv preprint arXiv:2202.05356},
  year={2022}
}

@article{chin2019regression,
	title={Regression adjustments for estimating the global treatment effect in experiments with interference},
	author={Chin, Alex},
	journal={Journal of Causal Inference},
	volume={7},
	number={2},
	year={2019},
	publisher={De Gruyter}
}

@article{SussmanAiroldi17,
	title={Elements of estimation theory for causal effects in the presence of network interference},
	author={Sussman, Daniel L and Airoldi, Edoardo M},
	journal={arXiv preprint arXiv:1702.03578},
	year={2017}
}

@article{BasseAiroldi15,
	title={Model-assisted design of experiments in the presence of network-correlated outcomes},
	author={Basse, Guillaume W and Airoldi, Edoardo M},
	journal={Biometrika},
	volume={105},
	number={4},
	pages={849--858},
	year={2018},
	publisher={Oxford University Press}
}

@article{cai2015social,
	title={Social networks and the decision to insure},
	author={Cai, Jing and De Janvry, Alain and Sadoulet, Elisabeth},
	journal={American Economic Journal: Applied Economics},
	volume={7},
	number={2},
	pages={81--108},
	year={2015}
}

@article{li2021causal,
	title={Causal Inference under Network Interference with Noise},
	author={Li, Wenrui and Sussman, Daniel L and Kolaczyk, Eric D},
	journal={arXiv preprint arXiv:2105.04518},
	year={2021}
}

@article{EcklesKarrerUgander17,
	title={Design and analysis of experiments in networks: Reducing bias from interference},
	author={Eckles, Dean and Karrer, Brian and Ugander, Johan},
	journal={Journal of Causal Inference},
	volume={5},
	number={1},
	year={2017},
	publisher={De Gruyter}
}

@inproceedings{GuiXuBhasinHan15,
	title={Network a/b testing: From sampling to estimation},
	author={Gui, Huan and Xu, Ya and Bhasin, Anmol and Han, Jiawei},
	booktitle={Proceedings of the 24th International Conference on World Wide Web},
	pages={399--409},
	year={2015},
	organization={International World Wide Web Conferences Steering Committee}
}

@article{savje2023causal,
  title={Causal inference with misspecified exposure mappings: separating definitions and assumptions},
  author={S{\"a}vje, Fredrik},
  journal={Biometrika},
  pages={asad019},
  year={2023},
  publisher={Oxford University Press}
}

@article{leung2022rate,
  title={Rate-optimal cluster-randomized designs for spatial interference},
  author={Leung, Michael P},
  journal={The Annals of Statistics},
  volume={50},
  number={5},
  pages={3064--3087},
  year={2022},
  publisher={Institute of Mathematical Statistics}
}

@article{hu2022switchback,
  title={Switchback experiments under geometric mixing},
  author={Hu, Yuchen and Wager, Stefan},
  journal={arXiv preprint arXiv:2209.00197},
  year={2022}
}

@article{horvitz1952generalization,
  title={A generalization of sampling without replacement from a finite universe},
  author={Horvitz, Daniel G and Thompson, Donovan J},
  journal={Journal of the American statistical Association},
  volume={47},
  number={260},
  pages={663--685},
  year={1952},
  publisher={Taylor \& Francis}
}

@article{bojinov2023design,
  title={Design and analysis of switchback experiments},
  author={Bojinov, Iavor and Simchi-Levi, David and Zhao, Jinglong},
  journal={Management Science},
  volume={69},
  number={7},
  pages={3759--3777},
  year={2023},
  publisher={INFORMS}
}

@inproceedings{blake2014marketplace,
  title={Why marketplace experimentation is harder than it seems: The role of test-control interference},
  author={Blake, Thomas and Coey, Dominic},
  booktitle={Proceedings of the fifteenth ACM conference on Economics and computation},
  pages={567--582},
  year={2014}
}

@article{glynn2020adaptive,
  title={Adaptive experimental design with temporal interference: A maximum likelihood approach},
  author={Glynn, Peter W and Johari, Ramesh and Rasouli, Mohammad},
  journal={Advances in Neural Information Processing Systems},
  volume={33},
  pages={15054--15064},
  year={2020}
}

@article{jagadeesan2020designs,
  title={Designs for estimating the treatment effect in networks with interference},
  author={Jagadeesan, Ravi and Pillai, Natesh S and Volfovsky, Alexander},
  year={2020}
}

@article{manski2013identification,
  title={Identification of treatment response with social interactions},
  author={Manski, Charles F},
  journal={The Econometrics Journal},
  volume={16},
  number={1},
  pages={S1--S23},
  year={2013},
  publisher={Oxford University Press Oxford, UK}
}

@inproceedings{toulis2013estimation,
  title={Estimation of causal peer influence effects},
  author={Toulis, Panos and Kao, Edward},
  booktitle={International conference on machine learning},
  pages={1489--1497},
  year={2013},
  organization={PMLR}
}

@article{forastiere2021identification,
  title={Identification and estimation of treatment and interference effects in observational studies on networks},
  author={Forastiere, Laura and Airoldi, Edoardo M and Mealli, Fabrizia},
  journal={Journal of the American Statistical Association},
  volume={116},
  number={534},
  pages={901--918},
  year={2021},
  publisher={Taylor \& Francis}
}

@article{basse2019randomization,
  title={Randomization tests of causal effects under interference},
  author={Basse, Guillaume W and Feller, Avi and Toulis, Panos},
  journal={Biometrika},
  volume={106},
  number={2},
  pages={487--494},
  year={2019},
  publisher={Oxford University Press}
}

@article{thomke2020building,
  title={Building a culture of experimentation},
  author={Thomke, Stefan},
  journal={Harvard Business Review},
  volume={98},
  number={2},
  pages={40--47},
  year={2020}
}

@article{kohavi2017surprising,
  title={The surprising power of online experiments},
  author={Kohavi, Ron and Thomke, Stefan},
  journal={Harvard business review},
  volume={95},
  number={5},
  pages={74--82},
  year={2017}
}

@article{ni2023design,
  title={Design of Panel Experiments with Spatial and Temporal Interference},
  author={Ni, Tu and Bojinov, Iavor and Zhao, Jinglong},
  journal={Available at SSRN 4466598},
  year={2023}
}

@inproceedings{jiang2016doubly,
  title={Doubly robust off-policy value evaluation for reinforcement learning},
  author={Jiang, Nan and Li, Lihong},
  booktitle={International Conference on Machine Learning},
  pages={652--661},
  year={2016},
  organization={PMLR}
}

@inproceedings{thomas2016data,
  title={Data-efficient off-policy policy evaluation for reinforcement learning},
  author={Thomas, Philip and Brunskill, Emma},
  booktitle={International Conference on Machine Learning},
  pages={2139--2148},
  year={2016},
  organization={PMLR}
}

@article{leung2023network,
  title={Network Cluster-Robust Inference},
  author={Leung, Michael P},
  journal={Econometrica},
  volume={91},
  number={2},
  pages={641--667},
  year={2023},
  publisher={Wiley Online Library}
}

@article{johari2022experimental,
  title={Experimental design in two-sided platforms: An analysis of bias},
  author={Johari, Ramesh and Li, Hannah and Liskovich, Inessa and Weintraub, Gabriel Y},
  journal={Management Science},
  volume={68},
  number={10},
  pages={7069--7089},
  year={2022},
  publisher={INFORMS}
}

@article{shi2023dynamic,
  title={Dynamic causal effects evaluation in a/b testing with a reinforcement learning framework},
  author={Shi, Chengchun and Wang, Xiaoyu and Luo, Shikai and Zhu, Hongtu and Ye, Jieping and Song, Rui},
  journal={Journal of the American Statistical Association},
  volume={118},
  number={543},
  pages={2059--2071},
  year={2023},
  publisher={Taylor \& Francis}
}

@inproceedings{ugander2013graph,
  title={Graph cluster randomization: Network exposure to multiple universes},
  author={Ugander, Johan and Karrer, Brian and Backstrom, Lars and Kleinberg, Jon},
  booktitle={Proceedings of the 19th ACM SIGKDD international conference on Knowledge discovery and data mining},
  pages={329--337},
  year={2013}
}

@article{farias2022markovian,
  title={Markovian interference in experiments},
  author={Farias, Vivek and Li, Andrew and Peng, Tianyi and Zheng, Andrew},
  journal={Advances in Neural Information Processing Systems},
  volume={35},
  pages={535--549},
  year={2022}
}

@article{larsen2023statistical,
  title={Statistical challenges in online controlled experiments: A review of a/b testing methodology},
  author={Larsen, Nicholas and Stallrich, Jonathan and Sengupta, Srijan and Deng, Alex and Kohavi, Ron and Stevens, Nathaniel T},
  journal={The American Statistician},
  pages={1--15},
  year={2023},
  publisher={Taylor \& Francis}
}

@article{aronow2017estimating,
	title={Estimating average causal effects under general interference, with application to a social network experiment},
	author={Aronow, Peter M and Samii, Cyrus and others},
	journal={The Annals of Applied Statistics},
	volume={11},
	number={4},
	pages={1912--1947},
	year={2017},
	publisher={Institute of Mathematical Statistics}
}
